\documentclass[onefignum,onetabnum]{siamart190516}

\usepackage[utf8]{inputenc} 
\usepackage[T1]{fontenc} 
\usepackage{amsmath,latexsym} 
\usepackage{mathtools}
\mathtoolsset{showonlyrefs=false,mathic,centercolon} 
\usepackage{amssymb,amsfonts} 
\usepackage{booktabs} 
\usepackage{enumitem}
\usepackage{babel}
\usepackage{algpseudocode} 
\usepackage{algorithm} 
\usepackage{bm} 
\usepackage{url}
\usepackage{tikz}
\usepackage{standalone} 
\usepackage{epstopdf}
\ifpdf
  \DeclareGraphicsExtensions{.eps,.pdf,.png,.jpg}
\else
  \DeclareGraphicsExtensions{.eps}
\fi

\usepackage[labelfont=sl,textfont=sl]{subcaption} 
\usepackage[font={small,sl},labelsep=period]{caption} 

\usetikzlibrary{shapes}
\usetikzlibrary{arrows.meta}

\newsiamremark{remark}{Remark}
\newsiamremark{hypothesis}{Hypothesis}
\crefname{hypothesis}{Hypothesis}{Hypotheses}
\newsiamthm{claim}{Claim}
\newsiamthm{assumption}{Assumption}

\headers{Model-Based DFO for Convex-Constrained Optimization}{M. Hough and L. Roberts}

\title{Model-Based Derivative-Free Methods for Convex-Constrained Optimization\thanks{Submitted to the editors 3rd December 2021.}}

\author{Matthew Hough\thanks{Department of Combinatorics and Optimization, University of Waterloo, 200 University Avenue W., Waterloo, ON N2L 3G1, Canada (\email{mhough@uwaterloo.ca})}
\and Lindon Roberts\thanks{Mathematical Sciences Institute, Building 145, Science Road, Australian National University, Canberra ACT 2601, Australia (\email{lindon.roberts@anu.edu.au}).}}

\newcommand{\secref}[1]{Section~\ref{#1}}
\newcommand{\defref}[1]{Definition~\ref{#1}}
\newcommand{\thmref}[1]{Theorem~\ref{#1}}
\newcommand{\lemref}[1]{Lemma~\ref{#1}}
\newcommand{\corref}[1]{Corollary~\ref{#1}}
\renewcommand{\algref}[1]{Algorithm~\ref{#1}}
\newcommand{\assref}[1]{Assumption~\ref{#1}}
\newcommand{\figref}[1]{Figure~\ref{#1}}

\newcommand{\remref}[1]{Remark~\ref{#1}}

\newcommand{\R}{\mathbb{R}} 
\newcommand{\N}{\mathbb{N}} 
\newcommand{\E}{\mathbb{E}} 
\newcommand{\bigO}{\mathcal{O}} 
\newcommand{\defeq}{:=} 
\newcommand{\grad}{\nabla} 
\newcommand{\proj}{\operatorname{proj}} 
\newcommand{\spanop}{\operatorname{span}} 
\newcommand{\norm}[1]{\left\lVert #1 \right\rVert} 
\newcommand{\abs}[1]{\left\lvert #1 \right\rvert} 
\newcommand{\bmat}[1]{\begin{bmatrix}#1\end{bmatrix}} 

\newcommand{\kappaef}{\kappa_{\textnormal{ef}}}
\newcommand{\kappaeg}{\kappa_{\textnormal{eg}}}
\newcommand{\gammainc}{\gamma_{\textnormal{inc}}}
\newcommand{\gammadec}{\gamma_{\textnormal{dec}}}
\newcommand{\flow}{f_{\textnormal{low}}}
\newcommand{\epsilonmc}{\epsilon_{\textnormal{mc}}}

\def\be{\begin{align}}
\def\ee{\end{align}}

\renewcommand{\b}[1]{\bm{#1}} 
\renewcommand{\t}[1]{\widetilde{#1}} 

\newcommand{\So}{\mathcal{S}} 
\newcommand{\Ps}{\mathcal{P}} 

\newcommand{\bzero}{\b{0}}
\newcommand{\bx}{\b{x}}
\newcommand{\bd}{\b{d}}

\newcommand{\by}{\b{y}}

\newcommand{\bs}{\b{s}}

\newcommand{\br}{\b{r}}
\newcommand{\bee}{\b{e}}
\newcommand{\bg}{\b{g}}
\newcommand{\bem}{\b{m}}

\algrenewcommand\algorithmicrequire{\textbf{Input:}}
\algrenewcommand\algorithmicensure{\textbf{Output:}}

\newcommand{\revision}[1]{{#1}}  

\ifpdf
\hypersetup{
  pdftitle={Model-Based Derivative-Free Methods for Convex-Constrained Optimization},
  pdfauthor={M. Hough, L. Roberts}
}
\fi

\begin{document}

\maketitle

\begin{abstract}
  We present a model-based derivative-free method for optimization subject to general convex constraints, which we assume are unrelaxable and accessed only through a projection operator that is cheap to evaluate.
  We prove global convergence and a worst-case complexity of $\bigO(\epsilon^{-2})$ iterations and objective evaluations for nonconvex functions, matching results for the unconstrained case.
  We introduce new, weaker requirements on model accuracy compared to existing methods.
  As a result, sufficiently accurate interpolation models can be constructed only using feasible points.
  We develop a comprehensive theory of interpolation set management in this regime for linear and composite linear models.
  We implement our approach for nonlinear least-squares problems and demonstrate strong practical performance compared to general-purpose solvers.
\end{abstract}






\section{Introduction}
The \revision{minimization of} a function when derivative information is not available, derivative-free optimization (DFO), is an area of growing research attention with many applications \cite{Conn2009,Audet2017,Larson2019}. 
The two most common approaches to DFO for finding local minima are direct search, where the objective is sampled in several directions to find descent, and model-based, which uses sampled function values to construct a local model for the objective within (typically) a trust-region framework.
The survey \cite{Larson2019} provides a detailed overview of these approaches.

Because of their relative simplicity, direct search methods have been extended to a wide variety of problem classes, including general nonsmooth \cite{Audet2006} or discontinuous \cite{Vicente2012} objectives, and those with discrete variables \cite{Audet2019}.
By contrast, the analysis of model-based DFO beyond general unconstrained optimization is less developed, largely due to the difficulty of forming appropriate notions of model accuracy.
Most commonly the local models are formed by interpolation, and model accuracy is guaranteed by requiring the interpolation set to satisfy particular geometric conditions.

In this work we consider model-based DFO for nonlinear optimization in the presence of convex constraints.
That is, we aim to solve
\begin{align}
  \min_{\bx\in C} f(\bx),  \label{eq_prob}
\end{align}
where $C \subseteq \R^n$ is closed and convex with nonempty interior.
We assume that the objective $f:\R^n\to\R$ is smooth and nonconvex, although access to its derivatives is not possible.
In our approach, the constraint set is made available to the algorithm only through a projection operator which we assume is cheap to compute.
To do this, we follow the approach from derivative-based trust-region methods \cite[Chapter 12]{Conn2000} in which the trust-region subproblem is solved using projected gradient descent.
We provide a first-order convergence analysis and worst-case complexity bounds for our approach, motivated by the analysis for derivative-based cubic regularization methods from \cite{Cartis2012b}.

Importantly, our algorithm is strictly feasible, and so $f$ is only ever evaluated at points in $C$.
This requires the development of new notions of model accuracy suitable when the interpolation set is restricted to lie in $C$.
To that end, we introduce a novel, comprehensive theory in which we generalize the notions of fully linear models and $\Lambda$-poised interpolation sets \cite{Conn2007,Conn2009a} to arbitrary convex sets.
This generalization allows us to overcome a known limitation of model-based DFO, namely that for some problems ``\textit{it may be impossible to obtain a fully linear model using only feasible points}'' \cite[p.~362]{Larson2019}.
Our theory is suitable for models based on linear interpolation, and so are particularly useful for composite minimization.
To that end, we provide numerical experiments on constrained nonlinear least-squares problems, where we show strong practical performance.

\subsection{Existing Works}
Trust-region methods for solving \eqref{eq_prob} in a derivative-based setting are well-studied, with \cite[Chapter 12]{Conn2000} providing a comprehensive development of measures of stationarity, the projected gradient method for solving the trust-region problem, and global convergence analysis.
A worst-case complexity analysis of this technique, where global convergence is guaranteed by adaptive cubic regularization was given in \cite{Cartis2012b}.
A related earlier work is \cite{Conn1993}, which considered the derivative-based trust-region method in the case of inexact gradients $\bg_k$ satisfying $\|\bg_k - \grad f(\bx_k)\| \leq \kappaeg \Delta_k$ at all iterations $k$ (c.f.~\eqref{eq_orig_fl}).
More recently, \cite{Bellavia2019a} considers high-order adaptive regularization methods for general constrained optimization problems subject to inexact objective and derivative information.

For model-based DFO, the works \cite{Conn2007,Conn2009a} introduced the relevant concepts of model accuracy (full linearity) and interpolation set geometry ($\Lambda$-poisedness) used to analyze unconstrained problems, and worse-case complexity analysis was developed (including for composite minimization) in \cite{Grapiglia2016,Garmanjani2016}.
A variety of model-based DFO approaches have been proposed for constrained optimization problems.
These vary by whether the constraints are known exactly or are black-box, and whether they are relaxable or not.
A comprehensive survey of these methods, categorized using the taxonomy of \cite{LeDigabel2015}, can be found in \cite[Section 7]{Larson2019}.
However, we particularly note model-based DFO methods for problems with unrelaxable constraints, including \cite{Powell2009,Gratton2011} and \cite[Section 6.3]{Wild2009} for bound constraints, and \cite{Gumma2014} for linear inequality constraints.
All of these works focus on the development of efficient algorithms and implementations, and do not have convergence analysis.
We extend these approaches to a significantly broader class of constraints by assuming a convex constraint set accessed only via projections, while providing theoretical guarantees of convergence.

The most similar existing work to ours is \cite{Conejo2013}, which introduces a model-based DFO method for convex constrained optimization.
Similar to our setting, access to the constraints is assumed to be through a projection operator and global convergence is analyzed through the framework of \cite[Chapter 12]{Conn2000} provided the local model is always fully linear in the sense of \cite{Conn2009a}.
Our work extends this by generalizing full linearity to general convex constraint sets, providing interpolation set management routines for ensuring these conditions, and having an algorithmic framework where models are (occasionally) allowed to not be fully linear.
We also provide a worst-case complexity analysis of our method (as well as global convergence) and numerical results for nonlinear least-squares problems.
We also note that the same first-order criticality measure was also used to study probabilistic direct search methods subject to linear constraints \cite{Gratton2019a}.

We will specialize our approach to composite minimization, such as nonlinear least-squares problems.
In the derivative-based setting, \cite{Porcelli2013} analyzes a Gauss-Newton method for nonlinear least-squares problems with bound constraints and inexact trust-region subproblem solutions.
In the DFO context, \cite{Zhang2010,DFOGN} analyze unconstrained model-based DFO methods for nonlinear least-squares problems, and \cite{Morini2017} studies a derivative-free linesearch method for nonlinear systems with convex constraints.

\subsection{Contributions}
The contributions of this paper are threefold.

First, we introduce a model-based DFO algorithm for \eqref{eq_prob}, naturally extending the standard framework of unconstrained model-based DFO \cite{Conn2009} to convex constraints, based on the approach in \cite{Conn2000}.
Our method only accesses constraint information through a projection operator, which is assumed to be cheap to evaluate.
Compared to \cite{Conejo2013}, our approach uses the same broad trust-region framework but with a different criticality measure \eqref{eq_opt_measure}.
We also introduce a new notion of fully linear models suitable for the constrained setting, while also permitting iterations where the local model is not accurate.
Our new fully linear definition (\defref{def_fl}) is similar to that of \cite{Conn2009a}, but is specifically adapted to the setting of convex constraints.
As in \cite{Conejo2013}, we prove global convergence to first-order critical points, but we augment this by showing a worst-case complexity bound of $\bigO(\epsilon^{-2})$ iterations to drive the criticality measure below $\epsilon$.
In the unconstrained case, our bound recovers the result of \cite{Garmanjani2016}, including the same dependency on problem dimension in the case of linear interpolation models.
Our approach allows a wider range of constraint sets than  \cite{Powell2009,Wild2009,Gratton2011,Gumma2014}, which also do not provide convergence analysis.

Second, in order to construct interpolation models satisfying our new notion of full linearity, we extend the theory of $\Lambda$-poisedness from \cite{Conn2007} to the constrained setting.
Our weaker notion of full linearity means that we can construct $\Lambda$-poised interpolation sets by only sampling from feasible points (overcoming the aforementioned limitation described in \cite[p.~362]{Larson2019}).
We show that linear interpolation based on $\Lambda$-poised sets produce fully linear models (in our new sense of both terms).
We also give concrete procedures to generate $\Lambda$-poised sets and validate whether or not a given set is $\Lambda$-poised.
Although linear interpolation models are sufficient for producing fully linear models, they have gained particular popularity in composite minimization (e.g.~\cite{Grapiglia2016,Garmanjani2016,DFOGN}), and so we show that $\Lambda$-poised sets with linear interpolation produces fully linear models for composite minimization.

Lastly, and motivated by the case of composite minimization, we use our approach to extend the state-of-the-art solver DFO-LS \cite{DFOLS} for nonlinear least-squares problems to handle general convex constraints via projection operators.\footnote{We note that DFO-LS could previously already handle box constraints using techniques from \cite{Powell2009}, but did not have any convergence guarantees.}
Testing on low-dimensional test problems with a variety of linear and Euclidean ball inequality constraints, we show that the new DFO-LS can outperform the general-purpose solvers COBYLA \cite{Powell1994}, which requires the functional form for the constraints, and NOMAD \cite{nomad2011}, which handles constraints via an extreme barrier.
These solvers differ based on how they access the constraints, and are not designed to exploit the least-squares structure of the objective.

\paragraph{Code Availability}
The latest version of DFO-LS with these improvements is available on Github.\footnote{Version 1.3, \url{https://github.com/numericalalgorithmsgroup/dfols}.}

\paragraph{Notation}
Throughout, we let $\proj_C$ denote the projection operator for the convex set $C\subseteq\R^n$ in the Euclidean norm.\footnote{This is a well-defined operator from $\R^n$ to $C$, as per \cite[Theorem 6.25]{Beck2017}, for instance.}
We use $\|\cdot\|$ to be the 2-norm for vectors and matrices, and define $B(\bx,\Delta) := \{\by\in\R^n : \|\by-\bx\| \leq \Delta\}$ to be the closed ball centered at $\bx\in\R^n$ of radius $\Delta\geq 0$.

\section{Algorithm Statement}
In this section, we outline the general model-based DFO algorithm for \eqref{eq_prob}.
Our framework closely mimics the unconstrained setting given in \cite[Chapter 10]{Conn2009}, with modifications for the constraints similar to the ones in the derivative-based trust-region setting \cite[Chapter 12]{Conn2000}.
We first outline our notion of first-order criticality for \eqref{eq_prob}, and then give our algorithm.

Throughout this work, we assume the following about our feasible set $C$.

\begin{assumption} \label{ass_feasible_set}
  The feasible set $C\subseteq\R^n$ is a closed, convex set with nonempty interior, $\operatorname{int} C \neq \emptyset$.
\end{assumption}

\subsection{Criticality Measure}
We will seek a first-order optimal solution to \eqref{eq_prob}.
To measure this, we use the following first-order criticality measure from \cite[Section 12.1.4]{Conn2000}, defined for all $\bx\in C$:
\begin{align}
    \pi^f(\bx) \defeq \abs{\min_{\substack{\bx+\bd \in C\\ \norm{\bd}\leq 1}} \grad f(\bx)^T \bd}. \label{eq_opt_measure} 
\end{align}
We will say that $\bx^*$ is a first-order critical point for \eqref{eq_prob} if $\pi^f(\bx^*)=0$ (see \cite[Theorem 12.1.6]{Conn2000}).
If $C=\R^n$, then we recover the first-order criticality measure $\pi^f(\bx)=\|\grad f(\bx)\|$.
As shown in \cite[Theorem 12.1.4]{Conn2000}, the minimizer of \eqref{eq_opt_measure} is of the form $\bd^* = \bx(t)-\bx$ for some $t \geq 0$ sufficiently large, where $\bx(t)$ is the projected gradient path, $\bx(t)\defeq \proj_C[\bx-t \grad f(\bx)]$.

However, in a DFO setting, we cannot measure $\pi^f(\bx)$ as it depends on $\grad f(\bx)$.
Instead, we will only have access to some approximate gradient $\bg\approx \grad f(\bx)$.
As a consequence, we will need the following result, which quantifies the sensitivity of the criticality measure \eqref{eq_opt_measure} to errors in the gradient $\grad f(\bx)$.
This is a modification of \cite[Theorem 3.4]{Cartis2012b}, which proves that $\pi^f(\bx)$ is Lipschitz continuous provided $\grad f(\bx)$ is also.

\begin{lemma} \label{lem_crit_perturb}
  Fix $\bx\in C$ and $\bg_1,\bg_2\in\R^n$, and define
  \begin{align}
      \pi_1(\bx) \defeq \abs{\min_{\substack{\bx+\bd \in C\\ \norm{\bd}\leq 1}} \bg_1^T \bd}, \qquad \text{and} \qquad \pi_2(\bx) \defeq \abs{\min_{\substack{\bx+\bd \in C\\ \norm{\bd}\leq 1}} \bg_2^T \bd}.
  \end{align}
  Then we have
  \begin{align}
      \abs{\pi_1(\bx) - \pi_2(\bx)} \leq \max_{\substack{\bx+\bd \in C\\ \norm{\bd}\leq 1}} \abs{(\bg_1-\bg_2)^T \bd}. \label{eq_crit_perturb}
  \end{align}
\end{lemma}
\begin{proof}
  Since $\bd=\bzero$ is feasible for both $\pi_1(\bx)$ and $\pi_2(\bx)$, the quantity inside the absolute values is always non-positive. 
  Hence, defining $\bd_i$ to be the minimizer corresponding to $\pi_i(\bx)$, we have $\pi_i(\bx) = -\bg_i^T \bd_i$.
  
  First, suppose that $\pi_1(\bx) \geq \pi_2(\bx)$.
  Then,
  \begin{align}
      \abs{\pi_1(\bx) - \pi_2(\bx)} = \pi_1(\bx) - \pi_2(\bx) &= -\bg_1^T \bd_1 + \bg_2^T \bd_2, \\
      &= -(\bg_1-\bg_2)^T \bd_1 + (\bg_2^T \bd_2 - \bg_2^T \bd_1), \\
      &\leq \abs{(\bg_1-\bg_2)^T \bd_1}, \label{eq_crit_tmp1}
  \end{align}
  where to get the last inequality we use $\bg_2^T \bd_2 \leq \bg_2^T \bd_1$ (since $\bd_2$ is the minimizer for $\pi_2(\bx)$ and $\bd_1$ is a feasible point for the same problem).
  
  Instead, if $\pi_1(\bx) < \pi_2(\bx)$, then by similar reasoning we have
  \begin{align}
      \abs{\pi_1(\bx) - \pi_2(\bx)} &\leq \abs{(\bg_2-\bg_1)^T \bd_2}. \label{eq_crit_tmp2}
  \end{align}
  Combining \eqref{eq_crit_tmp1} and \eqref{eq_crit_tmp2} proves \eqref{eq_crit_perturb}.
\end{proof}

We also note that applying Cauchy-Schwarz and $\|\bd\|\leq 1$ to \eqref{eq_crit_perturb} gives us the Lipschitz-type property
\begin{align}
    \abs{\pi_1(\bx) - \pi_2(\bx)} \leq \|\bg_1 - \bg_2\|. \label{eq_crit_lipschitz}
\end{align}

\subsection{Main Algorithm}
Our main algorithm follows a trust-region framework.
At each iteration $k$, we construct a local quadratic model which approximates the objective near the current iterate $\bx_k$:
\begin{align}
    f(\by) \approx m_k(\by) \defeq c_k + \bg_k^T (\by-\bx_k) + \frac{1}{2}(\by-\bx_k)^T H_k (\by-\bx_k), \label{eq_model}
\end{align}
where we hope this approximation is accurate when $\by\approx\bx_k$.
Following \cite[Chapter 12]{Conn2000}, we calculate a tentative new iterate $\bx_k+\bs_k$, where the step $\bs_k$ is an approximate minimizer of the trust-region subproblem
\begin{align}
    \min_{\substack{\bx+\bs \in C\\ \norm{\bs}\leq \Delta_k}} m_k(\bx_k+\bs), \label{eq_trs}
\end{align}
where the trust-region radius $\Delta_k > 0$ is updated at each iteration.
We accept this step (i.e.~set $\bx_{k+1}=\bx_k+\bs_k$) and increase $\Delta_k$ if this point sufficiently decreases $f$, otherwise we set $\bx_{k+1}=\bx_k$ and (possibly) decrease $\Delta_k$.

In the DFO setting, there is the possibility that the model \eqref{eq_model} is not a good approximation for the objective.
The usual notion of $m_k$ being a `good approximation' is the notion of a `fully linear' model (e.g.~\cite[Definition 6.1]{Conn2009}).
For the constrained case, we introduce a similar notion of full linearity, in part motivated by \eqref{eq_crit_perturb}.

\begin{definition}  \label{def_fl}
  The model $m_k$ \eqref{eq_model} is \emph{fully linear} in $B(\bx_k,\Delta_k)$ if there exist $\kappaef,\kappaeg>0$, independent of $k$, such that
  \begin{subequations} \label{eq_fl}
  \begin{align}
    \max_{\substack{\revision{\bx_k}+\bd \in C\\ \norm{\bd}\leq \Delta_k}} \abs{f(\bx_k+\bd) - m_k(\bx_k+\bd)} &\leq \kappaef \Delta_k^2, \label{eq_fl_f} \\
    \max_{\substack{\revision{\bx_k}+\bd \in C\\ \norm{\bd}\leq 1}} \abs{(\grad f(\bx_k)-\bg_k)^T \bd} &\leq \kappaeg \Delta_k. \label{eq_fl_g}
  \end{align}
  \end{subequations}
\end{definition}

In particular, we note that we require $\|\bd\|\leq\Delta_k$ in \eqref{eq_fl_f} but $\|\bd\|\leq 1$ in \eqref{eq_fl_g}.
In the unconstrained case $C=\R^n$, \defref{def_fl} is slightly weaker than the usual version of full linearity because it requires that $\grad m(\by) \approx \grad f(\by)$ only at $\by=\bx_k$, rather than all $\by\in B(\bx_k,\Delta_k)$.
We assume the existence of procedures which can verify whether $m_k$ is fully linear, and (if not) to make $m_k$ fully linear.
We describe such procedures in \secref{sec_models}.

The local model $m_k$ \eqref{eq_model} induces an approximate first-order criticality measure, namely
\begin{align}
    \pi^m(\bx) \defeq \abs{\min_{\substack{\bx+\bd \in C\\ \norm{\bd}\leq 1}} \bg_k^T \bd}. \label{eq_pik}
\end{align}
Our full linearity definition allows us to compare our approximate criticality measure $\pi^m_k \defeq \pi^m(\bx_k)$ with the true value $\pi^f_k \defeq \pi^f(\bx_k)$.

\begin{lemma} \label{lem_pi_comparison}
  If $m_k$ is fully linear in $B(\bx_k,\Delta_k)$, then the approximate criticality measure $\pi^m_k$ \eqref{eq_pik} satisfies 
  \begin{align}
      \abs{\pi^f_k - \pi^m_k} \leq \kappaeg \Delta_k. \label{eq_pi_comparison}
  \end{align}
\end{lemma}
\begin{proof}
  This immediately follows from \eqref{eq_fl_g} and \lemref{lem_crit_perturb}.
\end{proof}

In the unconstrained case, \eqref{eq_pi_comparison} is exactly  the fully linear condition \eqref{eq_fl_g}, and so there is nothing to show.

Our main algorithm, CDFO-TR, for solving \eqref{eq_prob} is given in \algref{alg_cdfotr}.
It is essentially the same as \cite[Algorithm 10.1]{Conn2009}, but we modify the criticality step to avoid having an inner loop, similar to \cite[Algorithm 3.1]{Garmanjani2016} and \cite[Algorithm 1]{DFBGN}, and replace references to the criticality measure $\|\bg_k\|$ with $\pi^m_k$.

\begin{algorithm}[tb]
\begin{algorithmic}[1]
\Require Starting point $\bx_0\in\R^n$ and trust-region radius $\Delta_0>0$.
\vspace{0.2em}
\Statex \underline{Parameters:} maximum trust-region radius $\Delta_{\max} \geq \Delta_0$, scaling factors $0 < \gammadec < 1 < \gammainc$, criticality constants $\epsilon_C,\mu>0$, and acceptance threshold $\eta\in(0,1)$.
\State Build an initial model $m_0$ \eqref{eq_model}.
\For{$k=0,1,2,\ldots$}
    \If{$\pi^m_k<\epsilon_C$ \textbf{and} \textit{($\pi^m_k<\mu^{-1}\Delta_k$ \textnormal{\textbf{or}} $m_k$ is not fully linear in $B(\bx_k,\Delta_k)$)}} \label{ln_main_start}
        \State \underline{Criticality step:} Set $\bx_{k+1}=\bx_k$. If $m_k$ is fully linear in $B(\bx_k,\Delta_k)$, set $\Delta_{k+1}=\gammadec\Delta_k$, otherwise set $\Delta_{k+1}=\Delta_k$. Construct $m_{k+1}$ to be fully linear in $B(\bx_{k+1},\Delta_{k+1})$.
    \Else \quad $\leftarrow$ \textit{$\pi^m_k\geq\epsilon_C$ or ($\pi^m_k\geq\mu^{-1}\Delta_k$ and $m_k$ is fully linear in $B(\bx_k,\Delta_k)$)}
        \State Approximately solve \eqref{eq_trs} to get a step $\bs_k$.
        \State Evaluate $f(\bx_k+\bs_k)$ and calculate ratio 
        \begin{align}
            \rho_k \defeq \frac{f(\bx_k)-f(\bx_k+\bs_k)}{m_k(\bx_k)-m_k(\bx_k+\bs_k)}. \label{eq_ratio_generic}
        \end{align}
        \If{$\rho_k \geq \eta$}
            \State \underline{Successful step:} Set $\bx_{k+1}=\bx_k+\bs_k$ and $\Delta_{k+1}=\min(\gammainc\Delta_k,\Delta_{\max})$. Form $m_{k+1}$ in any manner.
        \ElsIf{$m_k$ is not fully linear in $B(\bx_k,\Delta_k$)}
            \State \underline{Model-improving step:} Set $\bx_{k+1}=\bx_k$ and $\Delta_{k+1}=\Delta_k$, and construct $m_{k+1}$ to be fully linear in $B(\bx_{k+1},\Delta_{k+1})$.
        \Else \quad $\leftarrow$ \textit{$\rho_k<\eta$ and $m_k$ is fully linear in $B(\bx_k,\Delta_k)$}
            \State \underline{Unsuccessful step:} Set $\bx_{k+1}=\bx_k$ and $\Delta_{k+1}=\gammadec\Delta_k$. Form $m_{k+1}$ in any manner.
        \EndIf
    \EndIf \label{ln_main_end}
\EndFor
\end{algorithmic}
\caption{CDFO-TR: model-based DFO method for \eqref{eq_prob}.}
\label{alg_cdfotr}
\end{algorithm}

\section{Convergence \& Worst-Case Complexity} \label{sec_wcc}
We now prove convergence and provide a worst-case complexity analysis for \algref{alg_cdfotr}.
We will require the following (standard) assumptions.


\begin{assumption} \label{ass_smoothness}
  The objective function $f$ is bounded below by $\flow$ continuously differentiable.
  Furthermore, the gradient $\grad f$ is Lipschitz continuous with constant $L_{\grad f}$ in $\cup_k B(\bx_k,\Delta_{\max})$.
\end{assumption}


\begin{assumption} \label{ass_boundedhess}
  There exists $\kappa_H \geq 1$ such that $\|H_k\| \leq \kappa_H - 1$ for all $k$.
\end{assumption}


Lastly, we need an assumption about the accuracy with which the trust-region subproblem \eqref{eq_trs} is solved.

\begin{assumption} \label{ass_cdec}
  There exists a constant $c_1\in(0,1)$ such that the computed step $\bs_k$ satisfies $\bx_k+\bs_k\in C$, $\|\bs_k\|\leq\Delta_k$ and the generalized Cauchy decrease condition:
  \begin{align}
    m_k(\bx_k) - m_k(\bx_k + \bs_k) \geq c_1\pi^m_k \min\left(\frac{\pi^m_k}{1 + \|H_k\|}, \Delta_k, 1\right). \label{eq_2cdec}
  \end{align}
\end{assumption}

\assref{ass_cdec} is a slight modification to the unconstrained Cauchy decrease assumption \cite[Assumption AA.1]{Conn2000}.
We note that \assref{ass_cdec} is achievable using a Goldstein-type linesearch method \cite[Algorithm 12.2.2 \& Theorem 12.2.2]{Conn2000}.

\subsection{Convergence of \algref{alg_cdfotr}}
\begin{lemma} \label{lem_neverdecoupled}
  Suppose \assref{ass_smoothness} holds and the criticality step is not called in iteration $k$. 
  Then $\pi^m_k \geq \min(\epsilon_C,\Delta_k/\mu)$. 
  Also, if $\pi^f_k \geq \epsilon > 0$ then
  \begin{align}
    \pi^m_k \geq \epsilonmc \defeq \min\left(\epsilon_C, \frac{\epsilon}{1+\kappaeg\mu}\right) > 0. \label{eq_2pimkbnd}
  \end{align}
\end{lemma}
\begin{proof}
    This proof is based on \cite[Lemma 2.10]{DFBGN}.
    We first note that $\pi^m_k \geq \min(\epsilon_C,\Delta_k/\mu)$ is simply a necessary condition for not entering the criticality step.
    
    For the second part, suppose that the criticality step is not called in iteration $k$, and $\pi^f_k \geq \epsilon$ and $\pi^m_k < \epsilon_C$.
    Then $\Delta_k \leq \mu \pi^m_k$ and $m_k$ is fully linear in $B(\bx_k,\Delta_k)$, and so \lemref{lem_pi_comparison} gives
    \begin{align}
        \epsilon \leq \pi^f_k \leq \abs{\pi^f_k - \pi^m_k} + \pi^m_k \leq \kappaeg \Delta_k + \pi^m_k \leq (\kappaeg \mu+1)\pi^m_k,
    \end{align}
    and so $\pi^m_k \geq \epsilon / (1 + \kappaeg \mu)$.


\end{proof}

\begin{lemma} \label{lem_successfuliters}
  Suppose Assumptions~\ref{ass_smoothness}, \ref{ass_boundedhess} and \ref{ass_cdec} hold, $m_k$ is fully linear on $B(\bx_k,\Delta_k)$ and that
  \begin{align}
    \Delta_k \leq \min(c_0\pi^m_k, 1), \qquad \text{where} \qquad c_0 \defeq \min\left(\mu, \frac{1}{\kappa_H}, \frac{c_1(1-\eta)}{2\kappaef}\right). \label{eq_2dltaub}
  \end{align}
  Then the criticality step is not called in iteration $k$, and iteration $k$ is successful (i.e.~$\rho_k \geq \eta)$.
\end{lemma}

\begin{proof}
    This proof is similar to \cite[Lemma 2.6]{DFBGN}.
    Since $m_k$ is fully linear and $\Delta_k \leq c_0 \pi^m_k \leq \mu\pi^m_k$, the criticality step is not called.
    
    Then, Assumptions~\ref{ass_boundedhess} and \ref{ass_cdec} with $\Delta_k \leq \min(c_0 \pi^m_k, 1) \leq \min(\kappa_H^{-1} \pi^m_k, 1)$ imply
    \begin{align}
        m_k(\bx_k) - m_k(\bx_k+\bs_k) \geq c_1 \pi^m_k \min\left(\frac{\pi^m_k}{\kappa_H}, \Delta_k, 1\right) = c_1 \pi^m_k \Delta_k.
    \end{align}
    Now, since $m_k$ is fully linear, we apply \eqref{eq_fl_f} to get
    \begin{align}
        \abs{\rho_k - 1} &\leq \abs{\frac{f(\bx_k) - m_k(\bx_k)}{m_k(\bx_k) - m_k(\bx_k + \bs_k)}} + \abs{\frac{f(\bx_k + \bs_k) - m_k(\bx_k + \bs_k)}{m_k(\bx_k) - m_k(\bx_k + \bs_k)}}, \\
                     &\leq \frac{2\kappaef\Delta_k^2}{c_1\pi^m_k\Delta_k}, \\
                     &\leq 1-\eta, \label{eq_cauchy_tmp1}
    \end{align}
    where the last inequality follows from $\Delta_k \leq c_0 \pi^m_k \leq c_1 (1-\eta)\pi^m_k / (2\kappaef)$.
    That $\rho_k \geq \eta$ then follows from \eqref{eq_cauchy_tmp1}.

\end{proof}

\begin{lemma} \label{lem_deltalowerbound}
    Suppose that Assumptions~\ref{ass_smoothness}, \ref{ass_boundedhess} and \ref{ass_cdec} hold, and $\pi^f_k \geq \epsilon > 0$ for all $k\leq K$.
    Then $\Delta_k \geq \Delta_{\min}$ for all $k\leq K$, where
    \begin{align}
        \Delta_{\min} \defeq \min\left(\Delta_0, \frac{\gammadec \epsilon}{\kappaeg+\mu^{-1}}, \gammadec \mu \epsilonmc, \frac{\gammadec \epsilonmc}{\kappa_H}, \gammadec \left(\kappaeg + \frac{2\kappaef}{c_1 (1-\eta)}\right)^{-1} \epsilon, \gammadec\right). \label{eq_2deltalb}
    \end{align}
\end{lemma}
\begin{proof}
  This argument is similar to \cite[Lemma 3.8]{DFOGN}.
  To find a contradiction, suppose that $k'\leq K$ is the first iteration with $\Delta_{k'} < \Delta_{\min}$.
  Since $\Delta_{\min} \leq \Delta_0$, we must have $k'>0$.
  Given that $k'$ is the first such iteration, we have  $\Delta_{k'} < \Delta_{\min} \leq \Delta_{k'-1}$, and so iteration $(k'-1)$ must have either been a criticality step (with $m_{k'-1}$ fully linear) or an unsuccessful step; these are the only situations in which $\Delta_{k+1} < \Delta_k$. 
  
  First, suppose that iteration $(k'-1)$ had $m_{k'-1}$ fully linear in $B(\bx_{k'-1},\Delta_{k'-1})$ and the criticality step was called.
  By the entry condition for the criticality step, this means that $\pi^m_{k'-1} < \mu^{-1} \Delta_{k'-1}$.
  Combining these properties, we have
  \begin{align}
      \epsilon \leq \pi^f_{k'-1} \leq \abs{\pi^f_{k'-1} - \pi^m_{k'-1}} + \pi^m_{k'-1} \leq \kappaeg \Delta_{k'-1} + \mu^{-1}\Delta_{k'-1},
  \end{align}
  and so $\Delta_{k'-1} \geq \epsilon / (\kappaeg+\mu^{-1})$.
  Thus
  \begin{align}
      \Delta_{k'} = \gammadec \Delta_{k'-1} \geq \frac{\gammadec \epsilon}{\kappaeg+\mu^{-1}} \geq \Delta_{\min},
  \end{align}
  which is a contradiction.
  Thus iteration $(k'-1)$ must have been an unsuccessful step.
  
  For iteration $(k'-1)$ to have been an unsuccessful step, we must have had $\rho_{k'-1} < \eta$ and $m_{k'-1}$ fully linear in $B(\bx_{k'-1},\Delta_{k'-1})$.
  Also, since the criticality step is not called, we must have $\pi^m_{k'-1} \geq \epsilonmc$ by \lemref{lem_neverdecoupled}.
  
  
  Now, suppose that 
  \begin{align}
      \Delta_{k'-1} > \frac{c_1 (1-\eta)\pi^m_{k'-1}}{2\kappaef}. \label{eq_deltamin_tmp1}
  \end{align}
  Then by full linearity we have
  \begin{align}
      \epsilon \leq \kappaeg \Delta_{k'-1} + \pi^m_{k'-1} < \left(\kappaeg + \frac{2\kappaef}{c_1 (1-\eta)}\right)\Delta_{k'-1}. \label{eq_deltamin_tmp2}
  \end{align}
  However, since iteration $(k'-1)$ was unsuccessful, we have
  \begin{align}
      \Delta_{k'-1} = \gammadec^{-1} \Delta_{k'} < \gammadec^{-1} \Delta_{\min} \leq \left(\kappaeg + \frac{2\kappaef}{c_1 (1-\eta)}\right)^{-1} \epsilon,
  \end{align}
  contradicting \eqref{eq_deltamin_tmp2}, and hence \eqref{eq_deltamin_tmp1} must be false.
  
  Separately, we also have
  \begin{align}
      \Delta_{k'-1} < \gammadec^{-1} \Delta_{\min} \leq \min\left(\mu \epsilonmc, \frac{\epsilonmc}{\kappa_H}, 1\right),
  \end{align}
  and so altogether we have $\Delta_{k'-1} \leq \min(c_0 \pi^m_{k'-1}, 1)$.
  Hence by \lemref{lem_successfuliters} iteration $(k'-1)$ must have $\rho_{k'-1} \geq \eta$, a contradiction.
\end{proof}

We now proceed to show convergence of \algref{alg_cdfotr}.
The following results follow broadly the approach of \cite[Chapter 10]{Conn2009} and \cite{DFOGN}.

\begin{lemma} \label{lem_finite_successful_conv}
  Suppose \assref{ass_smoothness} holds and there are finitely many successful iterations. Then
  \begin{align}
  \lim_{k\to\infty} \pi^f_{k} = 0.
  \end{align}
\end{lemma}
\begin{proof}
  Let $k_0$ be the final successful iteration. 
  Then all iterations $k > k_0$ must either be criticality, model-improving or unsuccessful steps.
  If a model in some iteration $k > k_0$ is not fully linear, then \algref{alg_cdfotr} guarantees that the model in iteration $k+1$ is fully linear.
  Hence there are infinitely many iterations $k > k_0$ for which the model is fully linear.
  Such iterations (criticality or unsuccessful) always cause $\Delta_k$ to be reduced by a factor of $\gammadec$, and so $\lim_{k\to\infty} \Delta_k = 0$.

  Now suppose, for each $k > k_0$, that $i_k$ is the index of the first iteration after $k_0$ for which   $m_k$ is fully linear.
  Then by the reasoning above we must have $i_k\in\{k,k+1\}$.
  It follows that as $k \to \infty$,
  \begin{align}
    \norm{\bx_k - \bx_{i_k}} \leq \Delta_k \to 0. \label{eq_2xdiff}
  \end{align}
  Note that we may write
  \begin{equation}
    \pi^f_k \leq \abs{\pi^f_k - \pi^f_{i_k}} + \abs{\pi^f_{i_k} - \pi^m_{i_k}} + \pi^m_{i_k}.
  \end{equation}
  Since $\grad f$ is Lipschitz continuous, applying \cite[Theorem 3.4]{Cartis2012b} gives that $\abs{\pi^f_k - \pi^f_{i_k}} \leq M\norm{\bx_k - \bx_{i_k}}$ for some $M > 0$. 
  So from \eqref{eq_2xdiff}, the first term converges to zero.
  The second term converges to zero by using \lemref{lem_pi_comparison} to write $\abs{\pi^f_{i_k} - \pi^m_{i_k}} \leq \kappaeg\Delta_{i_k} \to 0$.
  Finally, the third term converges to zero, since if that were not the case, \lemref{lem_successfuliters}
  would contradict the fact that $k > k_0$ are not successful iterations. 
  It must then hold that $\pi^f_k\to 0$.
\end{proof}

\begin{lemma} \label{lem_deltatozero}
    Suppose that Assumptions~\ref{ass_smoothness}, \ref{ass_boundedhess} and \ref{ass_cdec} hold.
    Then
  \begin{align}
    \lim_{k\to\infty} \Delta_k = 0.
  \end{align}
\end{lemma}
\begin{proof}
  Let $\mathcal{S}$ be the set of successful iterations.
  Then the result when $\mathcal{S}$ is finite is shown in the proof of \lemref{lem_finite_successful_conv}.
  We will now consider the case when $\mathcal{S}$ is infinite. 
  
  For any $k \in \mathcal{S}$ we have from Assumptions~\ref{ass_boundedhess} and \ref{ass_cdec},
  \begin{align}
    f(\bx_k) - f(\bx_{k+1}) \geq \eta\left[m_k(\bx_k) - m_k(\bx_k + \bs_k)\right] \geq \eta c_1 \pi^m_k \min\left(\frac{\pi^m_k}{\kappa_H}, \Delta_k, 1\right).
  \end{align}
  Since the criticality step is not called for $k\in\mathcal{S}$, from \lemref{lem_neverdecoupled} we must have that $\pi^m_k \geq \min(\epsilon_C,\Delta_k/\mu)$; hence
  \begin{align}
    \sum_{k\in\mathcal{S}} f(\bx_k) - f(\bx_{k+1}) \geq \sum_{k\in\mathcal{S}} \eta c_1 \min(\epsilon_C,\Delta_k/\mu) \min\left(\frac{\min(\epsilon_C,\Delta_k/\mu)}{\kappa_H}, \Delta_k, 1\right) \geq 0.
  \end{align}
  By \assref{ass_smoothness}, the left-hand side of this inequality is bounded above by $f(\bx_0)-\flow$.
  Hence the right-hand side must be summable.
  Since $\mathcal{S}$ is infinite, this then implies that $\lim_{k\in\mathcal{S}} \Delta_k\to 0$.
  

  Now recall that the trust-region radius can be increased only during a successful iteration, and only by a ratio of at most $\gammainc$. Let $k \notin S$ be the index of any non-successful iteration after at least one successful iteration has occurred, and let $s_k \in\mathcal{S}$ be the last successful iteration before $k$. 
  Then $\Delta_k \leq \gammainc\Delta_{s_k}$. 
  Since from above we have that $\Delta_{s_k} \to 0$, then $\Delta_{k} \to 0$ for $k \notin S$ as $k \to \infty$, which proves the result in general.
\end{proof}

\begin{lemma} \label{lem_liminf}
Suppose that Assumptions~\ref{ass_smoothness}, \ref{ass_boundedhess} and \ref{ass_cdec} hold.
    Then
  \begin{align}
    \liminf_{k\to\infty} \pi^f_k = 0. \label{eq_liminf}
  \end{align}
\end{lemma}
\begin{proof}
    Let $\mathcal{S}$ be the set of successful iterations.
  Then the result when $\mathcal{S}$ is finite follows from \lemref{lem_finite_successful_conv}, so suppose $\mathcal{S}$ is infinite.
  
  Assume for contradiction that for all $k$ sufficiently large (e.g.~$k\geq k_0$),
  \begin{align}
    \pi^f_k \geq \epsilon, \label{eq_2pikdiverges}
  \end{align}
  for some $\epsilon > 0$. 
  From \lemref{lem_neverdecoupled}, $\pi^m_k \geq \epsilonmc > 0$ must also hold. 
  Now, consider some $k \in \mathcal{S}$.
  By Assumptions~\ref{ass_boundedhess} and \ref{ass_cdec}, \lemref{lem_deltalowerbound} and conditions \eqref{eq_2pikdiverges},
  \begin{align}
    f(\bx_k) - f(\bx_{k+1}) \geq \eta[m_k(\bx_k) - m_k(\bx_k + \bs_k)] \geq \eta c_1\epsilon_{mc}\min\left(\frac{\epsilonmc}{\kappa_H}, \Delta_{\min}, 1\right).
  \end{align}
  We can take the sum over all successful iterations from $k_0$ up to $k$:
  \begin{align}
    f(\bx_{k_0}) - f(\bx_{k+1}) = \sum_{\substack{j=k_0\\ j \in \mathcal{S}}}^k [f(\bx_j) - f(\bx_{j+1})] \geq \sigma_k \eta c_1\epsilonmc\left(\frac{\epsilonmc}{\kappa_H}, \Delta_{\min}, 1\right),
  \end{align}
  where $\sigma_k$ is the number of successful iterations between iterations $k_0$ and $k$. 
  However, since there are infinitely many successful iterations, we have that $\lim_{k\to\infty} \sigma_k = \infty$, and thus the difference between $f(\bx_{k_0})$ and $f(\bx_{k+1})$ is unbounded; so $f$ is unbounded below. 
  This contradicts \assref{ass_smoothness}, so we are done.
\end{proof}

\begin{theorem} \label{thm_limpif}
Suppose that Assumptions~\ref{ass_smoothness}, \ref{ass_boundedhess} and \ref{ass_cdec} hold.
    Then
  \begin{align}
    \lim_{k\to\infty} \pi^f_k = 0. \label{eq_2pifklim}
  \end{align}
  Moreover, $\lim\limits_{k\to\infty} \pi^m_k = 0$, cannot occur before \eqref{eq_2pifklim} holds.
\end{theorem}
\begin{proof}
  Assume for contradiction that there exists a subsequence of successful iterates indexed by $\{t_i\} \subseteq \mathcal{S}$ such that
  for some $\epsilon > 0$
  \begin{align}
    \pi^f_{t_i} \geq 2\epsilon > 0, \label{eq_2pikcontradiction}
  \end{align}
  for all $i$. 
  For each $t_i$, we know from \lemref{lem_liminf} that there exists a first successful iteration $\ell_i > t_i$ such that $\pi^f_{\ell_i} < \epsilon$. 
  Hence we have subsequences $\{t_i\}$ and $\{\ell_i\}$ of $\mathcal{S}$ such that such that
  \begin{align}
    \pi^f_{t_i} \geq 2\epsilon, \qquad \pi^f_k \geq \epsilon\ \text{ for $t_i < k < \ell_i$, \quad and \qquad $\pi^f_{\ell_i} < \epsilon$}. \label{eq_2pik1epsilon}
  \end{align}
  From \lemref{lem_neverdecoupled}, we must also have that $\pi^m_k \geq \epsilonmc > 0$ for $t_i \leq k < \ell_i$.
  Let us restrict our attention to the subsequence of successful iterations whose indices are in the set
  \begin{align}
  \mathcal{K} \defeq \bigcup_{i=1}^\infty\{k \in \mathcal{S} : t_i \leq k < \ell_i\},
  \end{align}
  where $t_i$ and $\ell_i$ belong to the two subsequences defined above. Since $\mathcal{K} \subseteq \mathcal{S}$, using $\pi^m_k \geq \epsilonmc$ and \assref{ass_boundedhess},
  for $k \in \mathcal{K}$ we have
  \begin{align}
    f(\bx_k) - f(\bx_{k+1}) \geq \eta[m_k(\bx_k) - m_k(\bx_k + \bs_k)] \geq \eta c_1\epsilonmc \min\left(\frac{\epsilonmc}{\kappa_H}, \Delta_k, 1\right) \geq 0. \label{eq_2fdiffbound}
  \end{align}
  \revision{From \lemref{lem_deltatozero} we have
  \begin{align}
    \label{eq_2limdelta}
    \lim_{\substack{k\to\infty\\k\in\mathcal{K}}} \Delta_k = 0.
  \end{align}}
  Applying \eqref{eq_2limdelta} to \eqref{eq_2fdiffbound}, we get that $\min(\epsilonmc/\kappa_H, \Delta_k, 1) = \Delta_k$ for $k$ large enough. 
  So,
  \begin{align}
    \label{eq_eq_2dltaub3}
    \Delta_k \leq \frac{1}{\eta c_1\epsilonmc}\left[f(\bx_k) - f(\bx_{k+1})\right],
  \end{align}
  for $k$ large enough. We can write
  \begin{align}
      \norm{\bx_{t_i} - \bx_{\ell_i}} \leq \sum_{\substack{j=t_i\\ j \in \mathcal{K}}}^{\ell_i - 1}\norm{\bx_j - \bx_{j+1}} \leq \sum_{\substack{j=t_i\\ j \in \mathcal{K}}}^{\ell_i - 1} \Delta_j &\leq \frac{1}{\eta c_1 \epsilonmc}\sum_{\substack{j=t_i\\ j\in\mathcal{K}}}^{\ell_i - 1}\left[f(\bx_j) - f(\bx_{j+1})\right], \\ 
      &= \frac{1}{\eta c_1 \epsilonmc}\left[f(\bx_{t_i}) - f(\bx_{\ell_i})\right], \label{eq_2successstepbound}
  \end{align}
observing that the first inequality comes from applying the triangle inequality to the
telescoping sum $\norm{\bx_{t_i} - \bx_{t_{i+1}} + \bx_{t_{i+1}} - \bx_{t_{i+2}} + \cdots}$, and the final inequality comes from \eqref{eq_eq_2dltaub3}. 
We can use the assumption that $f$ is bounded below with the fact that the sequence $\{f(\bx_k)\}$ is monotone decreasing, to say that $\{f(\bx_k)\}$ converges and is thus Cauchy. 
So the right hand side of \eqref{eq_2successstepbound} goes to zero as $i\to\infty$, and hence $\norm{\bx_{t_i} - \bx_{\ell_i}}$ also goes to zero as $i\to\infty$. 

Using \lemref{lem_crit_perturb}, we have that  $\abs{\pi^f_{t_i} - \pi^f_{\ell_i}} \leq M\norm{\bx_{t_i} - \bx_{\ell_i}} \to 0$.
But from the definitions of $\{t_i\}$ and $\{\ell_i\}$, this is a contradiction to $\abs{\pi^f_{t_i} - \pi^f_{\ell_i}} \geq \epsilon > 0$, a consequence of \eqref{eq_2pik1epsilon}, and we are done. 
\end{proof}


\subsection{Worst-Case Complexity}
We now analyze the worst-case complexity of \algref{alg_cdfotr}.
That is, we bound the number of iterations and objective evaluations until $\pi^f_k < \epsilon$ is first achieved, in terms of the desired optimality level $\epsilon$.
To that end, we define $i_\epsilon$ to be the last iteration before $\pi^f_k<\epsilon$ for the first time (which exists for all $\epsilon>0$ by \lemref{lem_liminf}).
That is, $\pi^f_k \geq \epsilon$ for all $k\leq i_{\epsilon}$ and $\pi^f_{i_\epsilon+1} < \epsilon$.



\begin{lemma} \label{lem_Sbound}
    Suppose that Assumptions~\ref{ass_smoothness}, \ref{ass_boundedhess} and \ref{ass_cdec} hold, and let $S_{i_\epsilon}$ be the set of iterations up to $i_\epsilon$ that are successful. Then
  \begin{align}
    \abs{S_{i_\epsilon}} \leq \frac{f(\bx_0)-\flow}{\eta c_1}\max\left(\epsilonmc^{-1}\Delta_{\min}^{-1}, \epsilonmc^{-1}\right), \label{eq_3Sbound}
  \end{align}
  where $\epsilonmc$ is defined in \eqref{eq_2pimkbnd} and $\Delta_{\min}$ in \eqref{eq_2deltalb}.
\end{lemma}
\begin{proof}
    
  For all $k \in S_{i_\epsilon}$, we have the sufficient decrease condition:
  \begin{align}
    f(\bx_k) - f(\bx_{k+1}) \geq \eta[m_k(\bx_k) - m_k(\bx_k + \bs_k)] \geq \eta c_1 \pi^m_k\min\left(\frac{\pi^m_k}{\kappa_H},\Delta_{k}, 1\right). \label{eq_complexity_tmp1}
  \end{align}
  From Lemmas~\ref{lem_neverdecoupled} and \ref{lem_deltalowerbound}, we have that $\pi^m_k \geq \epsilonmc$ and $\Delta_k \geq \Delta_{\min}$ for all $k\leq i_{\epsilon}$, and hence \eqref{eq_complexity_tmp1} becomes
  \begin{align}
    f(\bx_k) - f(\bx_{k+1}) \geq \eta c_1 \epsilonmc\min\left(\frac{\epsilonmc}{\kappa_H},\Delta_{\min}, 1\right), \qquad \forall k\in S_{i_{\epsilon}}. \label{eq_3sdec}
  \end{align}
  Summing \eqref{eq_3sdec} over all $k \in S_{i_\epsilon}$, and noting that $\flow \leq f(\bx_k) \leq f(\bx_0)$, we obtain
  \begin{align}
    f(\bx_0)-\flow \geq 
    \sum_{k\in S_{i_{\epsilon}}} f(\bx_k) - f(\bx_{k+1}) \geq \abs{S_{i_\epsilon}}\eta c_1 \epsilonmc\min\left(\frac{\epsilonmc}{\kappa_H}, \Delta_{\min}, 1\right).
  \end{align}
  Finally, \eqref{eq_3Sbound} follows from $\Delta_{\min} \leq \gammadec \epsilonmc/\kappa_H \leq \epsilonmc/\kappa_H$.
\end{proof}

It remains to count the number of iterations of \algref{alg_cdfotr} that are not successful. Counting until iteration $i_\epsilon$
(inclusive), we let
\begin{enumerate}
  \item $C_{i_\epsilon}^M$ be the set of criticality step iterations $k \leq i_\epsilon$ for which $\Delta_k$ is not reduced. This iteration occurs at most once per sequence of criticality steps;
  \item $C_{i_\epsilon}^U$ be the set of criticality step iterations $k \leq i_\epsilon$ for which $\Delta_k$ is reduced. This amounts to every criticality iteration except perhaps the single iteration where $\Delta_k$ is not reduced;
  \item $M_{i_\epsilon}$ be the set of iterations where the model-improvement step is called; and
  \item $U_{i_\epsilon}$ be the set of unsuccessful iterations.
\end{enumerate}

\begin{lemma} \label{lem_countbounds}
  Suppose that Assumptions~\ref{ass_smoothness}, \ref{ass_boundedhess} and \ref{ass_cdec} hold.
  We have the bounds
  \begin{align}
    \abs{U_{i_\epsilon}} + \abs{C_{i_\epsilon}^U} &\leq \frac{\log \gammainc}{\abs{\log\gammadec}}\abs{S_{i_\epsilon}} + \frac{\log(\Delta_0/\Delta_{\min})}{\abs{\log\gammadec}}, \label{eq_3cb1} \\
    \abs{C^M_{i_\epsilon}} + \abs{M_{i_\epsilon}} &\leq \abs{C^U_{i_{\epsilon}}} + \abs{S_{i_{\epsilon}}} + \abs{U_{i_{\epsilon}}} + 1. \label{eq_3cb2}
  \end{align}
\end{lemma}
\begin{proof}
  Note that if $k \in S_{i_\epsilon}$, we set $\Delta_{k+1} \leq \gammainc\Delta_k$; and if $k \in U_{i_\epsilon}$, we set $\Delta_{k+1} = \gammadec\Delta_k$.
  For iterations in $C_{i_\epsilon}^U$, we also set $\Delta_{k+1} = \gammadec\Delta_k$. Thus,
  \begin{align}
    \gammainc^{\abs{S_{i_\epsilon}}}\gammadec^{\abs{U_{i_\epsilon}} + \abs{C^U_{i_\epsilon}}} \geq \frac{\Delta_{i_{\epsilon}}}{\Delta_0} \geq \frac{\Delta_{\min}}{\Delta_0},
  \end{align}
  where the second inequality follows from \lemref{lem_deltalowerbound}.
  The first result \eqref{eq_3cb1} follows by taking the log of both sides.
  
  
  Now, if $k\in C^M_{i_\epsilon} \cup M_{i_{\epsilon}}$, then the model $m_k$ was not fully linear, but the mechanisms of \algref{alg_cdfotr} ensure that $m_{k+1}$ is fully linear.
  Hence iteration $k+1$ must be a criticality step in which $\Delta_{k+1}$ is reduced, successful or unsuccessful.
  Thus $k+1\in C^U_{i_{\epsilon}} \cup S_{i_{\epsilon}} \cup U_{i_{\epsilon}}$, or $k+1>i_{\epsilon}$ (i.e.~$k=i_{\epsilon}$), and \eqref{eq_3cb2} follows.
\end{proof}

To get an overall complexity bound, we will lastly need to assume that the algorithm parameter $\epsilon_C$ is not much smaller than the desired tolerance $\epsilon$.

\begin{assumption}
  \label{ass_complexity}
  The algorithm parameter $\epsilon_C \geq c_2 \epsilon$ for some constant $c_2 > 0$.
\end{assumption}
As noted in \cite{DFOGN}, Assumption~\ref{ass_complexity} can be easily satisfied by choosing $\epsilon_C$ appropriately in Algorithm~\ref{alg_cdfotr}.

\begin{theorem}
  \label{thm_complexity}
  Suppose that Assumptions~\ref{ass_smoothness}, \ref{ass_boundedhess}, \ref{ass_cdec} and \ref{ass_complexity} hold.
  The number of iterations $i_\epsilon$ until $\pi^f_{i_\epsilon + 1} < \epsilon$ is at most
  \begin{align}
      & \frac{f(\bx_0)-\flow}{\eta c_1}\left(2 + \frac{2\log \gammainc}{\abs{\log\gammadec}}\right)\max\left(\frac{1}{c_3 \epsilon \Delta_0}, \frac{1}{c_3 c_4 \epsilon^2}, \frac{1}{c_3 \gammadec \epsilon}, \frac{1}{c_3 \epsilon}\right) \nonumber \\
      &\qquad\qquad\qquad + \frac{2}{\abs{\log\gammadec}}\max\left(0, \log\left(\frac{\Delta_0}{c_4 \epsilon}\right), \log\left(\gammadec^{-1}\Delta_0\right)\right) + 1, \label{eqn_complexity-bound}
  \end{align}
  where $c_3 := \min(c_2,(1+\kappaeg\mu)^{-1})$ and
  \begin{align}
      c_4 := \gammadec \min\left(\frac{1}{\kappaeg+\mu^{-1}}, \mu c_3, \frac{c_3}{\kappa_H}, \left(\kappaeg + \frac{2\kappaef}{c_1 (1-\eta)}\right)^{-1}\right).
  \end{align}
 \end{theorem}
 \begin{proof}
    We first observe that, using \lemref{lem_countbounds}, the total number of iterations can be written as
    \begin{align}
        i_\epsilon &= \abs{S_{i_\epsilon}} + \abs{M_{i_\epsilon}} + \abs{U_{i_\epsilon}} + \abs{C^M_{i_\epsilon}} + \abs{C^U_{i_\epsilon}}, \\
        &\leq 2 \left(\abs{S_{i_\epsilon}} + \abs{U_{i_\epsilon}} + \abs{C^U_{i_\epsilon}}\right) + 1, \\
        &\leq \left(2 + \frac{2\log \gammainc}{\abs{\log\gammadec}}\right) \abs{S_{i_\epsilon}} + \frac{2\log(\Delta_0/\Delta_{\min})}{\abs{\log\gammadec}} + 1.
    \end{align}
    Hence, we apply Lemma~\ref{lem_Sbound} to obtain
    \begin{align}
        i_{\epsilon} \leq \left(2 + \frac{2\log \gammainc}{\abs{\log\gammadec}}\right) \frac{f(\bx_0)-\flow}{\eta c_1}\max\left(\epsilonmc^{-1}\Delta_{\min}^{-1}, \epsilonmc^{-1}\right) + \frac{2\log(\Delta_0/\Delta_{\min})}{\abs{\log\gammadec}} + 1. \label{cplxty_intermediate_s1}
    \end{align}
    From \lemref{lem_neverdecoupled} and \assref{ass_complexity}, we have $\epsilonmc = \min(\epsilon_C, \epsilon/(1+\kappaeg\mu)) \geq c_3 \epsilon$. 
    Similarly, \eqref{eq_2deltalb} gives
    \begin{align}
        \Delta_{\min} &\geq \min\left(\Delta_0, \frac{\gammadec \epsilon}{\kappaeg+\mu^{-1}}, \gammadec \mu c_3 \epsilon, \frac{\gammadec c_3 \epsilon}{\kappa_H}, \gammadec \left(\kappaeg + \frac{2\kappaef}{c_1 (1-\eta)}\right)^{-1} \epsilon, \gammadec\right), \\
        &= \min\left(\Delta_0, c_4 \epsilon, \gammadec\right),
    \end{align}
    and so
    \begin{align}
        \log\left(\frac{\Delta_0}{\Delta_{\min}}\right) \leq \max\left(0, \log\left(\frac{\Delta_0}{c_4  \epsilon}\right), \log\left(\gammadec^{-1}\Delta_0\right)\right), \label{cplxty_intermediate_s2}
    \end{align}
    and
    \begin{align}
        \max\left(\epsilonmc^{-1}\Delta_{\min}^{-1}, \epsilonmc^{-1}\right) &\leq \max\left(c_3^{-1} \epsilon^{-1} \Delta_{\min}^{-1}, c_3^{-1} \epsilon^{-1}\right), \\
        &\leq \max\left(\frac{1}{c_3 \epsilon \Delta_0}, \frac{1}{c_3 c_4 \epsilon^2}, \frac{1}{c_3 \gammadec \epsilon}, \frac{1}{c_3 \epsilon}\right). \label{cplxty_intermediate_s3}
    \end{align}
    As desired, \eqref{eqn_complexity-bound} follows from \eqref{cplxty_intermediate_s1}, \eqref{cplxty_intermediate_s2} and \eqref{cplxty_intermediate_s3}.
\end{proof}

\begin{corollary}
    \label{cor_complexity}
    Suppose the assumptions of \thmref{thm_complexity} hold.
    For $\epsilon \in (0,1]$, the number of iterations of \algref{alg_cdfotr} until $\pi^f_k < \epsilon$ for the first time is at most $\bigO(\kappa_H \kappa_d^2 \epsilon^{-2})$, where $\kappa_d := \max(\kappaef, \kappaeg)$. 
\end{corollary}
\begin{proof}
    We note that $c_0^{-1} = \bigO(\max(\kappa_H, \kappa_d))$ from \lemref{lem_successfuliters}, and $c_3^{-1} = \bigO(\kappa_d)$ and $c_4^{-1} = \bigO(\kappa_H \kappa_d)$ from \thmref{thm_complexity}.
    The result then follows immediately from \thmref{thm_complexity}.
\end{proof}

\begin{remark} \label{rem_complexity_bounds}
    We will see in \secref{sec_models} that if we are using linear or composite linear interpolation models, then we require at most $\bigO(n)$ objective evaluations per iteration.
    Further, we show in Theorems~\ref{thm_fl_new} and \ref{thm_fl_new_composite} that under standard smoothness conditions we can get $\kappa_d = \bigO(n)$ for linear interpolation and $\kappa_d=\bigO(n^2)$ for composite linear interpolation.
    We also have Hessian bound $\kappa_H=1$ for linear models (since there is no model Hessian) and the realistic estimate $\kappa_H = \bigO(\kappa_d) = \bigO(n^2)$ for composite linear interpolation (see \eqref{eq_composite_hess_bound}).
\end{remark}

Combining \corref{cor_complexity} with \remref{rem_complexity_bounds} gives a worst-case complexity bound of $\bigO(n^2 \epsilon^{-2})$ iterations and \revision{$\bigO(n^3 \epsilon^{-2})$} evaluations for linear interpolation models, and $\bigO(n^6 \epsilon^{-2})$ iterations and $\bigO(n^7 \epsilon^{-2})$ evaluations for composite linear interpolation.
These results match the first-order complexity bounds for unconstrained model-based DFO both for fully linear models \cite{Garmanjani2016} and nonlinear least-squares  with composite linear models \cite{DFOGN}.\footnote{Note that different authors use different conventions for estimating $\kappa_d$ in terms of $n$; see \cite[Remark 3.20]{DFOGN} for details.}

\section{Construction of Fully Linear Models} \label{sec_models}
In this section we consider the requirement for full linearity in the constrained case (\defref{def_fl}).
As above, we generalize the standard approach for constructing fully linear models, namely the concept of $\Lambda$-poisedness.
Specifically, in this section, we look at approaches for constructing linear interpolation models which satisfy the fully linear requirement.
We show that an approach based on maximizing Lagrange polynomials, similar to the unconstrained approach, allows us to construct fully linear models and/or verify if an existing model is fully linear.
We also extend this approach to linear interpolation for composite function minimization, for which linear models are sufficient to yield high-quality model-based DFO algorithms \cite{DFOGN}.

\subsection{Unconstrained Linear Interpolation}

Suppose we have a set of $n+1$ points in $\R^n$, denoted $\{\by_0,\ldots,\by_n\}$.
If we have evaluated our objective $f$ at all these points, we can construct a linear model (c.f.~\eqref{eq_model})
\begin{align}
    f(\by) \approx m(\by) \defeq c + \bg^T (\by-\by_0), \label{eq_lin_interp}
\end{align}
by enforcing the interpolation conditions
\begin{align}
    f(\by_t) = m(\by_t), \qquad \forall t=1,\ldots,n.
\end{align}
This reduces to the linear system
\begin{align}
    M \bmat{ c \\ \bg } \defeq \bmat{1 & (\by_0-\by_0)^T \\ \vdots & \vdots \\ 1 & (\by_n-\by_0)^T} \bmat{ c \\ \bg} = \bmat{f(\by_0) \\ \vdots \\ f(\by_n)}. \label{eq_interp_system}
\end{align}
This system is invertible provided the directions $\{\by_t-\by_0 : t=1,\ldots,n\}$ are linearly independent \cite[Section 2.3]{Conn2009}.

Given our interpolation set, the Lagrange polynomials associated with with this set are the linear functions
\begin{align}
    \ell_t(\by) \defeq c_t + \bg_t^T (\by-\by_0), \qquad \forall t=0,\ldots,n,
\end{align}
defined by $\ell_t(\by_s) = \delta_{s,t}$.
These conditions give
\begin{align}
    M \bmat{ c_t \\ \bg_t } = \bee_t, \qquad \forall t=0,\ldots,n,  \label{eq_lagrange_system}
\end{align}
where $\bee_t\in\R^{n+1}$ is the $t$-th coordinate vector.
More details of this construction are given in \cite[Section 3.2]{Conn2009}.

In the case of linear interpolation, the key result regarding full linearity is the following.

\begin{theorem} \label{thm_orig_fl}
  Suppose $f$ satisfies \assref{ass_smoothness}.
  Then if $\|\by_t-\by_0\| \leq \Delta$ for all $t=1,\ldots,n$ and there exists $\Lambda\geq 1$ such that
  \begin{align}
      \max_{t=0,\ldots,n} |\ell_t(\by)| \leq \Lambda, \qquad \forall \by\in B(\by_0,\Delta), \label{eq_orig_poised}
  \end{align}
  then the linear interpolation model \eqref{eq_lin_interp} satisfies
  \begin{subequations} \label{eq_orig_fl}
  \begin{align}
      \abs{f(\by) - m(\by)} &\leq \left(\frac{(3 + (n+1)\Lambda)L_{\grad f}}{2}\right)\Delta^2, \\
      \norm{\grad f(\by) - \grad m(\by)} &\leq \left(\frac{(2 + (n+1)\Lambda)L_{\grad f}}{2}\right)\Delta,
  \end{align}
  \end{subequations}
  for all $\by\in B(\by_0,\Delta)$.
\end{theorem}
\begin{proof}
  This comes from \cite[Theorems 2.11, 2.12 \& 3.14]{Conn2009}.
\end{proof}

The condition \eqref{eq_orig_poised} is called $\Lambda$-poisedness (e.g.~\cite[Definition 3.6]{Conn2009}) and \eqref{eq_orig_fl} gives the standard notion of full linearity (which immediately implies our fully linear notion---\defref{def_fl}---for any constraint set $C$) with specific values for $\kappaef$ and $\kappaeg$.
We also note that \eqref{eq_orig_poised} is closely related to the condition number of $M$ \cite[eq.~(3.11) \& Theorem 3.14]{Conn2009}, and that $\ell_0(\by_0)=1$ implies $\Lambda \geq 1$.

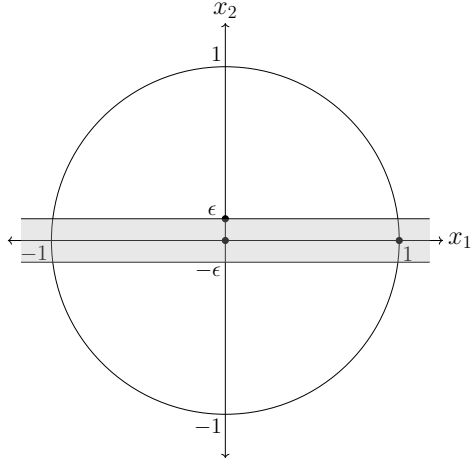
\begin{figure}[tb]
    \centering
    \begin{tikzpicture}

\draw[<->] (-5,0) -- (5,0) node[anchor=west] {\LARGE $x_1$};
\draw[<->] (0,-5) -- (0,5) node[anchor=south] {\LARGE $x_2$};

\node[] at (4.2,-0.3) {\Large $1$};
\node[] at (-4.4,-0.3) {\Large $-1$};
\node[] at (-0.2,4.3) {\Large $1$};
\node[] at (-0.4,-4.3) {\Large $-1$};

\draw (4,0) arc (0:360:4);

\filldraw (0,0) circle (2pt) node {};
\filldraw (4,0) circle (2pt) node {};
\filldraw (0,0.5) circle (2pt) node {};

\draw[] (-4.7,0.5) -- (4.7,0.5);
\draw[] (-4.7,-0.5) -- (4.7,-0.5);
\node[] at (-0.3,0.7) {\Large $\epsilon$};
\node[] at (-0.4,-0.7) {\Large $-\epsilon$};

\fill[gray!60,opacity=0.3] (-4.7,-0.5) -- (-4.7,0.5) -- (4.7,0.5) -- (4.7,-0.5) -- cycle;

\end{tikzpicture}
  \caption{Example feasible region and interpolation points with bad $\Lambda$-poisedness using \eqref{eq_orig_poised} but good $\Lambda$-poisedness using \eqref{eq_poised_demo}.}
  \label{fig_bad_poisedness}
\end{figure}

\paragraph{Constrained Example}
The standard notion of $\Lambda$-poisedness \eqref{eq_orig_poised} allows us to construct fully linear models (\defref{def_fl}) via linear interpolation.
However, to achieve this with a small value of $\Lambda$, it may be necessary to include interpolation points outside the feasible set $C$.
For example, consider the case $C=\{(x_1,x_2) : |x_2| \leq \epsilon\} \subset\R^2$ (for some $\epsilon\ll 1$) and the interpolation set $\{\bzero, (1, 0), (0, \epsilon)\}$ in $B(\bzero,1)$.
This is depicted in \figref{fig_bad_poisedness}.
In the ball $B(\bzero,1)$, these points are $\Lambda$-poised with $\Lambda \sim \epsilon^{-1}$, and so the corresponding fully linear constants \eqref{eq_orig_fl} will be similarly large.
However, if we instead only consider the size of the Lagrange polynomials only within the feasible set,
\begin{align}
    \max_{t=0,\ldots,n} \abs{\ell_t(\by)} \leq \Lambda, \qquad \forall \by\in C \cap B(\by_0,\Delta), \label{eq_poised_demo}
\end{align}
then it is not hard\footnote{This can be shown by relaxing $B(\by_0,\Delta)\cap C$ in \eqref{eq_poised_demo} to $\{(y_1,y_2) : |y_1|\leq 1, \, |y_2|\leq\epsilon\}$.} to show that our interpolation set satisfies \eqref{eq_poised_demo} with $\Lambda \leq 3$ for any choice of $\epsilon$. 

This example demonstrates that, in the constrained case, we are likely to get improved error bounds if we consider the size of the Lagrange polynomials not over the whole ball $B(\by_0,\Delta)$, but only the feasible subset of this ball.

\subsection{$\Lambda$-poisedness in the constrained case}
We now prove a generalization of \thmref{thm_orig_fl} suitable for problems with convex constraints.
To do this, we first introduce a notion of $\Lambda$-poisedness suitable for this setting.

\begin{definition} \label{def_new_poised}
  Given $\Lambda\geq 1$, an interpolation set $\{\by_0,\ldots,\by_n\} \subset C$ is $\Lambda$-poised for linear interpolation in $B(\by_0,\Delta)\cap C$ if $\{\by_t-\by_0 : t=1,\ldots,n\}$ is linearly independent and \begin{align}
      \max_{t=0,\ldots,n} \abs{\ell_t(\by)} \leq \Lambda, \qquad \forall \by\in C \cap B(\by_0,\min(\Delta,1)). \label{eq_new_poised}
  \end{align}
\end{definition}

Note that this is slightly weaker than \eqref{eq_poised_demo} as we only have $\by\in B(\by_0,\min(\Delta,1))$ rather than $B(\by_0,\Delta)$ (c.f.~\cite[Lemma 10.25]{Conn2009}).

\begin{lemma} \label{lem_fl_intermediate}
  Suppose $f$ satisfies \assref{ass_smoothness} and $C$ satisfies \assref{ass_feasible_set}.
  Then if $\{\by_0,\ldots,\by_n\}$ is $\Lambda$-poised for linear interpolation in $B(\by_0,\Delta)\cap C$ and $\|\by_t-\by_0\| \leq \beta \min(\Delta, 1)$ for all $t=1,\ldots,n$ and some $\beta>0$, we have
  \begin{align}
    \abs{(\by-\by_0)^T (\bg-\grad f(\by_0))} \leq \frac{n \Lambda L_{\grad f} \beta^2}{2}\min(\Delta,1)^2, \label{eq_fl_intermediate}
  \end{align}
  for all $\by\in B(\by_0, \min(\Delta,1))\cap C$.
\end{lemma}
\begin{proof}
  The first row of \eqref{eq_interp_system} gives $c=f(\by_0)$.
  Thus the remaining interpolation conditions give
  \begin{align}
      (\by_t-\by_0)^T \bg = f(\by_t) - f(\by_0), \qquad \forall t=1,\ldots,n.
  \end{align}
  Therefore, for all $t=1,\ldots,n$ we have
  \begin{align}
      \abs{(\by_t-\by_0)^T (\bg-\grad f(\by_0))} &= \abs{f(\by_t) - f(\by_0) - \grad f(\by_0)^T(\by_t-\by_0)}, \\
      &\leq \frac{\revision{L_{\grad f}}}{2}\|\by_t-\by_0\|^2 \leq \frac{L_{\grad f} \beta^2}{2}\min(\Delta,1)^2, \label{eq_interp_tmp1}
  \end{align}
  which follows from \assref{ass_smoothness} and $\|\by_t-\by_0\|\leq\beta\min(\Delta,1)$.
  
  Now, fix any $\by\in B(\by_0,\min(\Delta,1))\cap C$.
  Since $\{\by_t-\by_0 : t=1,\ldots,n\}$ forms a basis of $\R^n$ (as it is linearly independent), there exist  $\alpha_1(\by),\ldots,\alpha_n(\by)\in\R$ such that
  \begin{align}
      \by - \by_0 = \sum_{t=1}^{n} \alpha_t(\by) (\by_t-\by_0). \label{eq_interp_tmp2}
  \end{align}
  Hence from \eqref{eq_interp_tmp1} we get
  \begin{align}
      \abs{(\by-\by_0)^T (\bg-\grad f(\by_0))} &\leq \sum_{t=1}^{n}\abs{\alpha_t(\by)(\by_t-\by_0)^T (\bg-\grad f(\by_0))}, \\
      &\leq \frac{L_{\grad f} \beta^2}{2} \left(\sum_{t=1}^{n} \abs{\alpha_t(\by)}\right)\min(\Delta,1)^2. \label{eq_interp_tmp4}
  \end{align}
  We now construct a bound on each $|\alpha_t(\by)|$ in terms of $\Lambda$.
  For $t=1,\ldots,n$, the condition $\ell_t(\by_0)=0$ gives $c_t=0$, and so the remaining interpolation conditions $\ell_t(\by_s)=\delta_{s,t}$ give
  \begin{align}
      L\bg_t \defeq \bmat{(\by_1-\by_0)^T \\ \vdots \\ (\by_n-\by_0)^T} \bg_t = \bee_t,
  \end{align}
  where here $\bee_t\in\R^n$ is the $t$-th coordinate vector.
  Hence the $\Lambda$-poisedness condition gives
  \begin{align}
      \abs{\ell_t(\by)} = \abs{\bg_t^T (\by-\by_0)} = \abs{\bee_t^T L^{-T}(\by-\by_0)} \leq \Lambda. \label{eq_interp_tmp3}
  \end{align}
  Separately, we can find the values $\alpha_t(\by)$ by solving the linear system associated with \eqref{eq_interp_tmp2}, which is
  \begin{align}
      \underbrace{\bmat{(\by_1-\by_0) & \cdots & (\by_n-\by_0)}}_{=L^T} \bmat{\alpha_1(\by) \\ \vdots \\ \alpha_n(\by)} = \by-\by_0.
  \end{align}
  Combining this with \eqref{eq_interp_tmp3} we get
  \begin{align}
      \abs{\alpha_t(\by)} = \abs{[L^{-T}(\by-\by_0)]_t} = \abs{\bee_t^T L^{-T}(\by-\by_0)} \leq \Lambda.
  \end{align}
  This and \eqref{eq_interp_tmp4} gives us the result.
\end{proof}

As we might hope, this notion of $\Lambda$-poisedness guarantees that our interpolation models are fully linear.

\begin{theorem} \label{thm_fl_new}
  Suppose $f$ and $C$ satisfy Assumptions~\ref{ass_smoothness} and \ref{ass_feasible_set} respectively.
  Then if $\{\by_0,\ldots,\by_n\}$ is $\Lambda$-poised for linear interpolation in $B(\by_0,\Delta)\cap C$ and $\|\by_t-\by_0\| \leq \beta \min(\Delta, 1)$ for all $t=1,\ldots,n$ and some $\beta>0$, then the interpolation model \eqref{eq_lin_interp} is fully linear in $B(\by_0,\Delta)$ (in the sense of \defref{def_fl}) with constants
  \begin{align}
      \kappaef = \frac{L_{\grad f}(n\Lambda \beta^2 + 1)}{2}, \qquad \text{and} \qquad \kappaeg = \frac{n\Lambda L_{\grad f} \beta^2}{2},
  \end{align}
  in \eqref{eq_fl}.
\end{theorem}
\begin{proof}
    Fix $\by\in B(\by_0,\Delta)\cap C$.
    We first derive $\kappaef$ by considering the cases $\Delta\geq 1$ and $\Delta<1$ separately.
    If $\Delta\geq 1$ then $\hat{\by} := \by_0 + \Delta^{-1}(\by-\by_0) \in B(\by_0,1)=B(\by_0,\min(\Delta,1))$ and $\hat{\by}\in C$ since $C$ is convex.
    We can then apply \lemref{lem_fl_intermediate} to get
    \begin{align}
        \abs{(\by-\by_0)^T (\bg-\grad f(\by_0))} &= \Delta \abs{(\hat{\by}-\by_0)^T (\bg-\grad f(\by_0))}, \\
        &\leq \frac{n \Lambda L_{\grad f} \beta^2}{2}\Delta \min(\Delta,1)^2, \\
        &\leq \frac{n \Lambda L_{\grad f} \beta^2}{2} \Delta^2, \label{eq_poised_tmp1}
    \end{align}
    where the last inequality follows from  $\min(\Delta,1)^2  = 1 \leq \Delta$.
    Instead, if $\Delta<1$ then $\by\in B(\by_0,\min(\Delta,1))\cap C$ and so we apply \lemref{lem_fl_intermediate} directly, giving
    \begin{align}
        \abs{(\by-\by_0)^T (\bg-\grad f(\by_0))} \leq \frac{n \Lambda L_{\grad f} \beta^2}{2} \min(\Delta,1)^2 = \frac{n \Lambda L_{\grad f} \beta^2}{2} \Delta^2. \label{eq_poised_tmp2}
    \end{align}
    Finally, regardless of the size of $\Delta$, we again use $c=f(\by_0)$ (as in the proof of \lemref{lem_fl_intermediate}) with either \eqref{eq_poised_tmp1} or \eqref{eq_poised_tmp2} to get
  \begin{align}
      \abs{f(\by) - m(\by)} &= \abs{f(\by) - f(\by_0) - \bg^T(\by-\by_0)}, \\
      &\leq \abs{f(\by) - f(\by_0) - \grad f(\by_0)^T(\by-\by_0)} + \abs{(\by-\by_0)^T (\bg-\grad f(\by_0))}, \\
      &\leq \frac{L_{\grad f}}{2}\Delta^2 + \frac{n\Lambda L_{\grad f} \beta^2}{2} \Delta^2,
  \end{align}
  and we get the desired value of $\kappaef$.
  
  To get $\kappaeg$ we now fix an arbitrary $\t{\by}\in B(\by_0,1)\cap C$ and again consider the cases $\Delta\geq 1$ and $\Delta<1$ separately.
  First, if $\Delta\geq 1$, then $\t{\by} \in B(\by_0,\min(\Delta,1))\cap C$, and applying \lemref{lem_fl_intermediate} we get
  \begin{align}
      \abs{(\t{\by}-\by_0)^T (\bg-\grad f(\by_0))} \leq \frac{n \Lambda L_{\grad f} \beta^2}{2}\min(\Delta,1)^2 \leq \frac{n \Lambda L_{\grad f} \beta^2}{2} \Delta, \label{eq_interp_tmp6a}
  \end{align}
  since $\min(\Delta,1)^2 = 1 \leq \Delta$.
  Alternatively, if $\Delta<1$ then the convexity of $C$ implies that $\hat{\by}\defeq \by_0+\Delta(\t{\by}-\by_0)\in B(\by_0,\Delta)\cap C = B(\by_0,\min(\Delta,1))\cap C$.
  Again we apply \lemref{lem_fl_intermediate} and get
  \begin{align}
      \abs{(\t{\by}-\by_0)^T (\bg-\grad f(\by_0))} &= \Delta^{-1}\abs{(\hat{\by}-\by_0)^T (\bg-\grad f(\by_0))}, \\
      &\leq \frac{n \Lambda L_{\grad f} \beta^2}{2}\Delta^{-1}\min(\Delta,1)^2 = \frac{n \Lambda L_{\grad f} \beta^2}{2}\Delta. \label{eq_interp_tmp6b}
  \end{align}
  The value for $\kappaeg$ then follows from \eqref{eq_interp_tmp6a} and \eqref{eq_interp_tmp6b}.
\end{proof}

\begin{remark}
    The proof of \thmref{thm_orig_fl} strongly relies on the fact that $\|M\|=\|M^T\|$ (in the first line of \cite[Theorem 3.14]{Conn2009}), where $M$ arises from the interpolation model and $M^T$ arises from the construction of Lagrange polynomials.
    However, a similar statement is not necessarily true when we consider the matrix operator norm restricted to inputs in $C$, and so to prove \thmref{thm_fl_new} we introduce a novel approach to achieve very similar results.
\end{remark}

\thmref{thm_fl_new} gives exactly the type of result we want: that in order to achieve accurate models inside the feasible region, we only need to control the size of $\ell_t(\by)$ inside the feasible region.
As our example in \figref{fig_bad_poisedness} illustrates, this can be a substantially easier requirement to satisfy, particularly with the restriction that all interpolation points must be feasible.

With regards to our worst-case complexity bound \corref{cor_complexity}, we note that \thmref{thm_fl_new} gives $\kappaef,\kappaeg = \bigO(n L_{\grad f} \Lambda)$, which is the same size as in \thmref{thm_orig_fl} (and so in terms of $n$, \revision{$L_{\grad f}$} and $\Lambda$ our results match the complexity bounds from \cite{Garmanjani2016}).

\subsection{Constructing $\Lambda$-Poised Sets}
We now turn our attention to procedures which check whether or not a set is $\Lambda$-poised (in the sense of \defref{def_new_poised}), and to modify a non-poised set to ensure that it is $\Lambda$-poised.
For simplicity of notation, in this section we define $r\defeq \beta\min(\Delta,1)$.

\subsubsection{Verifying $\Lambda$-poisedness} \label{subsec_verify}
Verifying whether or not a set $\{\by_0,\ldots,\by_n\}\subset C$ is $\Lambda$-poised per \defref{def_new_poised} is straightforward, based on the procedure:
\begin{enumerate}
    \item Verify if the set $\{\by_t-\by_0 : t=1,\ldots,n\}$ is linearly independent (e.g.~via QR factorization, or checking if the matrix $M$ in \eqref{eq_lagrange_system} is invertible);
    \item Construct each Lagrange polynomial $\ell_t(\by)$ for $t=0,\ldots,n$ by solving the corresponding system \eqref{eq_lagrange_system};
    \item For each $t=0,\ldots,n$, calculate
    \begin{align}
        \max_{\by\in C\cap B(\by_0,\min(\Delta,1))} |\ell_t(\by)|,
    \end{align}
    and verify if this quantity is at most $\Lambda$.
    This reduces to solving two smooth convex optimization problems---minimizing $\pm\ell_t(\by)$---for each $t$ (e.g.~using the projected gradient method \cite[Theorem 12.1.4]{Conn2000}), at least to sufficient accuracy that it is known whether the optimal value is at most $\Lambda$.
\end{enumerate}
Further, to guarantee fully linear models (using \thmref{thm_fl_new}) we also need to ensure each interpolation point $\by_t\in C\cap B(\by_0,r)$: clearly $\|\by_t-\by_0\| \leq r$ is easily checked, and $\by_t\in C$ can be checked by verifying $\proj_C(\by_t)=\by_t$.

\subsubsection{Constructing Affinely Independent Points}
If the process to verify $\Lambda$-poisedness outlined in \secref{subsec_verify} fails, the first reason for this would be that the interpolation set $\{\by_0,\ldots,\by_n\}$ does not produce $n$ linearly independent directions $\{\by_t-\by_0 : t=1,\ldots,n\}$ or is not contained in $C\cap B(\by_0,r)$.
In this case, we need to replace our interpolation set by one satisfying these properties.
At this stage we are not concerned with verifying the $\Lambda$-poisedness property \eqref{eq_new_poised}; this is addressed in \secref{subsec_geometry}.

\begin{figure}[tb]
  \centering
  \begin{subfigure}[b]{0.48\textwidth}
    \begin{tikzpicture}

\draw[<->] (-5,0) -- (5,0) node[anchor=west] {\LARGE $x_1$};
\draw[<->] (0,-5) -- (0,5) node[anchor=south] {\LARGE $x_2$};

\node[] at (3.6,0.3) {\Large $\bm{y_0}+r\bm{e}_1$};
\node[] at (1.1,4) {\Large $\bm{y_0}+r\bm{e}_2$};

\node[] at (-0.4,0.3) {\Large $\bm{y_0}$};

\draw [dashed,-{Latex[length=2mm]}] (0,4) to [out=300,in=60] (0.1,0.1);
\draw [dashed,-{Latex[length=2mm]}] (4,0) to [out=210,in=330] (0.1,-0.1);
\node[] at (1.5,2) {\Large $\operatorname{proj}_C$};
\node[] at (2,-1) {\Large $\operatorname{proj}_C$};

\filldraw (0,0) circle (2pt) node {};
\filldraw (4,0) circle (2pt) node {};
\filldraw (0,4) circle (2pt) node {};

\node[] at (-2, -2) {\Large $C$};
\fill[gray!60,opacity=0.3] (-4.7,0) -- (-4.7,-4.7) -- (0,-4.7) -- (0,0) -- cycle;

\end{tikzpicture}
    \caption{$\proj_C(\by_0+r\bee_t)=\by_0$ for all $t=1,\ldots,n$.}
  \end{subfigure}
  ~
  \begin{subfigure}[b]{0.48\textwidth}
    \begin{tikzpicture}

\draw[<->] (-5,0) -- (5,0) node[anchor=west] {\LARGE $x_1$};
\draw[<->] (0,-5) -- (0,5) node[anchor=south] {\LARGE $x_2$};

\node[] at (3.6,0.3) {\Large $\bm{y_0}+r\bm{e}_1$};
\node[] at (1.1,4) {\Large $\bm{y_0}+r\bm{e}_2$};

\node[] at (0.4,0.3) {\Large $\bm{y_0}$};

\draw [dashed,-{Latex[length=2mm]}] (0,4) to (-1.9,2.1);
\draw [dashed,-{Latex[length=2mm]}] (4,0) to (2.1,-1.9);
\node[] at (-1.5,3.3) {\Large $\operatorname{proj}_C$};
\node[] at (3.3,-1.5) {\Large $\operatorname{proj}_C$};

\filldraw (0,0) circle (2pt) node {};
\filldraw (4,0) circle (2pt) node {};
\filldraw (0,4) circle (2pt) node {};

\filldraw (-2,2) circle (2pt) node {};
\filldraw (2,-2) circle (2pt) node {};

\node[] at (-2, -2) {\Large $C$};
\draw[] (-4.7,4.7) to (4.7,-4.7);
\fill[gray!60,opacity=0.3] (-4.7,4.7) -- (4.7,-4.7) -- (-4.7,-4.7) -- cycle;

\end{tikzpicture}
    \caption{$\{\by_0, \proj_C(\by_0+r\bee_1), \proj_C(\by_0+r\bee_2)\}$ colinear.}
  \end{subfigure}
  \caption{Examples of potential difficulties constructing affinely independent points using $\proj_C(\by_0+r\bee_t)$ for $t=1,\ldots,n$.}
  \label{fig_bad_initial_sets}
\end{figure}
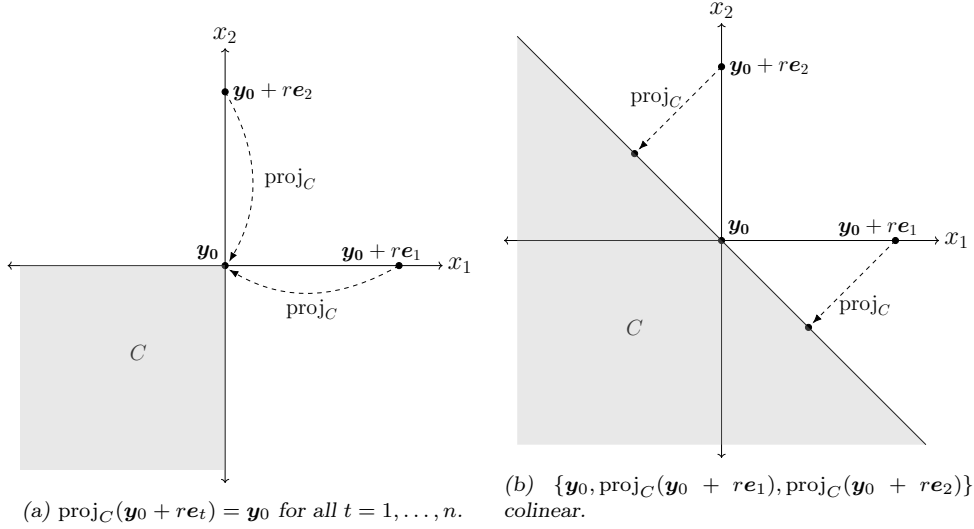

In the unconstrained case, it is easy to create a suitable replacement set: simply taking $\by_0$ and $\by_t=\by_0+r\bee_t$ for $t=1,\ldots,n$ is sufficient.
However, this is more difficult when $C\neq \R^n$, as two undesirable situations may arise, even from sampling directions orthogonal to $\by_0$:
\begin{itemize}
    \item $\proj_C(\by_0+r\bee_t)$ might be equal to $\by_0$ for some or all $t=1,\ldots,n$. This arises, for example if $C=\{\bx\in\R^2 : x_1,x_2 \leq 0\}$ and $\by_0=(0,0)$.
    \item $\proj_C(\by_0+r\bee_{t_1})-\by_0$ and $\proj_C(\by_0+r\bee_{t_2})-\by_0$ may not be linearly independent even if $t_1\neq t_2$. For example, if $C = \{\bx\in\R^2 : x_1 + x_2 \leq 0\}$ and $\by_0=(0,0)$, then $\proj_C(\by_0+r\bee_1)-\by_0 = (r/2,-r/2)$ is parallel to $\proj_C(\by_0+r\bee_2)-\by_0 = (-r/2, r/2)$.
\end{itemize}
These situations are depicted in \figref{fig_bad_initial_sets}. 

To address these difficulties and guarantee that we get a linearly independent set of directions after projecting onto $C$, we replace the perturbations $\{r\bee_t : t=1,\ldots,n\}$ with a countably infinite set of directions which are dense on the sphere $\{\bd\in\R^n : \|\bd\|=r\}$, and continue adding points of the form $\proj_C(\by_0+\bd)$ to the interpolation set until we reach $n+1$ points with $n$ linearly independent directions.
The precise algorithm is given in \algref{alg_dense_sampling}.
We note that there are many constructions for sampling dense directions on the unit sphere, such as the methods used to generate poll directions for Mesh Adaptive Direct Search (MADS) \cite{Audet2006,Abramson2009}.

\begin{algorithm}[tb]
\begin{algorithmic}[1]
\Require Constraint set $C\subset\R^n$ with projection operator $\proj_C:\R^n\to C$, starting point $\by_0\in C$, and sampling radius $r>0$.
\vspace{0.2em}
\State Construct a sequence $(\bd_k)_{k\in\N}$ of unit vectors $\bd_k\in\R^n$ which are dense on the unit sphere, and initialize $\mathcal{Y}=\{\by_0\}$ and $\mathcal{D}=\emptyset$.
\For{$k=1,2,\ldots$}
    \State Set $\by_k = \proj_C(\by_0+r\bd_k)$.
    \If{$\by_k-\by_0 \notin\spanop\mathcal{D}$}
        \State Set $\mathcal{Y}=\mathcal{Y}\cup\{\by_k\}$ and $\mathcal{D}=\mathcal{D}\cup\{\by_k-\by_0\}$.
        \If{$\mathcal{Y}$ has $n+1$ elements}
            \State \textbf{return} $\mathcal{Y}$.
        \EndIf
    \EndIf
\EndFor
\end{algorithmic}
\caption{Method for constructing affinely independent points.}
\label{alg_dense_sampling}
\end{algorithm}

\begin{theorem} \label{thm_dense_sampling}
  If $C$ satisfies \assref{ass_feasible_set} and $\by_0\in C$, then \algref{alg_dense_sampling} terminates in finite time, returning a set $\{\by_0,\ldots,\by_n\}\subset C\cap B(\by_0,r)$ for which $\{\by_t-\by_0 : t=1,\ldots,n\}$ is linearly independent.
\end{theorem}
\begin{proof}
    Provided that \algref{alg_dense_sampling} terminates, then by construction we know that it will return a set $\{\by_0,\ldots,\by_n\}\subset C$ for which $\{\by_t-\by_0:t=1,\ldots,n\}$ is linearly independent.
    Since $\by_0\in C$ we also have
    \begin{align}
        \|\by_k - \by_0\| = \|\proj_C(\by_0 + r\bd_k) - \proj_C(\by_0)\| \leq r \|\bd_k\| = r,
    \end{align}
    where the inequality follows from the non-expansiveness of the projection operator (e.g.~\cite[Theorem 6.42]{Beck2017}).
    Hence each $\by_t\in B(\by_0,r)$ and the result follows.
    
    It remains to prove that \algref{alg_dense_sampling} terminates in finite time.
    We do this by showing that, for any $\mathcal{D}$ with $|\mathcal{D}|<n$, there will be infinitely many iterations for which $\by_k-\by_0\notin\spanop\mathcal{D}$.
    On the first such iteration, $\mathcal{Y}$ and $\mathcal{D}$ gain one more element, and then either $|\mathcal{D}|=n$ and we are done, or we are again in the case $|\mathcal{D}|<n$ (and so again there are infinitely many future iterations in which $\mathcal{D}$ can be expanded).
    
    To this end, suppose at any iteration that $|\mathcal{D}|<n$.
    Since $C$ has nonempty interior and $\by_0\in C$, there exists $\t{\by}\in \operatorname{int}(C) \cap B(\by_0,r/2)$ with $\t{\by}\notin\by_0+\spanop\mathcal{D}$ (since $\by_0+\spanop\mathcal{D}$ is a set of measure zero).
    Let $\t{\by}_{\mathcal{D}}$ be the closest point to $\t{\by}$ in $\by_0+\spanop\mathcal{D}$, and so $\t{\by}_{\mathcal{D}} \neq \t{\by}$ and $ \t{\by}-\t{\by}_{\mathcal{D}} \in (\spanop\mathcal{D})^{\perp}$.
    The non-expansiveness of projections implies $\|\t{\by}_{\mathcal{D}}-\by_0\| \leq \|\t{\by}-\by_0\| \leq r/2$
    
    We will show that, since $\t{\by}\in C$, points near to the line between $\t{\by}_{\mathcal{D}}$ and $\t{\by}$, but which are closer to $\t{\by}$, give directions for which $\proj_C(\by_0+r\bd) \notin \by_0+\spanop\mathcal{D}$.
    The line between $\t{\by}_{\mathcal{D}}$ and $\t{\by}$ can be parametrized by $\alpha \mapsto \t{\by}_{\mathcal{D}} + \alpha(\t{\by}-\t{\by}_{\mathcal{D}})$.
    We can compute
    \begin{align}
        \|\t{\by}_{\mathcal{D}} + \alpha(\t{\by}-\t{\by}_{\mathcal{D}}) - \by_0\|^2 &= \alpha^2 \|\t{\by}-\t{\by}_{\mathcal{D}}\|^2 + 2\alpha (\t{\by}-\t{\by}_{\mathcal{D}})^T(\t{\by}_{\mathcal{D}}-\by_0) + \|\t{\by}_{\mathcal{D}}-\by_0\|^2.
    \end{align}
    Since $\|\t{\by}_{\mathcal{D}}-\by_0\|^2 \leq r^2/4 < r^2$ and $\|\t{\by}-\t{\by}_{\mathcal{D}}\|^2>0$, there exists a unique $\alpha^*>0$ such that $ \|\t{\by}_{\mathcal{D}} + \alpha(\t{\by}-\t{\by}_{\mathcal{D}}) - \by_0\|^2=r^2$.
    Let $\t{\by}^* \defeq \t{\by}_{\mathcal{D}} + \alpha^*(\t{\by}-\t{\by}_{\mathcal{D}})$ and $\bd^* \defeq \t{\by}^* - \by_0$.
    For this point, we have $\t{\by}^*-\t{\by}_{\mathcal{D}} = \alpha^*(\t{\by}-\t{\by}_{\mathcal{D}}) \in (\spanop\mathcal{D})^{\perp}$, and so $\t{\by}_{\mathcal{D}}=\proj_{\by_0+\spanop\mathcal{D}}(\t{\by}^*)$.
    
    In fact, we must have $\alpha^*>1$, since if  $\alpha\in[0,1]$ we have
    \begin{align}
        \|\t{\by}_{\mathcal{D}} + \alpha(\t{\by}-\t{\by}_{\mathcal{D}}) - \by_0\| \leq (1-\alpha)\|\t{\by}_{\mathcal{D}}-\by_0\| + \alpha\|\t{\by}-\by_0\| \leq r/2.
    \end{align}
    Thus
    \begin{align}
        \|\t{\by}^*-\t{\by}\| &= |1-\alpha^*|\, \|\t{\by}-\t{\by}_{\mathcal{D}}\| = \alpha^* \|\t{\by}-\t{\by}_{\mathcal{D}}\| - \|\t{\by}-\t{\by}_{\mathcal{D}}\|, \label{eq_poised_tmp4a} \\
        &<
        \alpha^* \|\t{\by}-\t{\by}_{\mathcal{D}}\| = \|\t{\by}^*-\t{\by}_{\mathcal{D}}\|. \label{eq_poised_tmp4}
    \end{align}
    Hence $\t{\by}^*$ is strictly closer to $\t{\by}\in C$ than $\t{\by}_{\mathcal{D}}$, and so we conclude that $\proj_C(\t{\by}^*) \notin \by_0+\spanop\mathcal{D}$.
    
    Lastly, consider any $\by\in B(\t{\by}^*,\|\t{\by}-\t{\by}_{\mathcal{D}}\|/4)$.
    Then we have
    \begin{align}
        \|\by-\proj_{\by_0+\spanop\mathcal{D}}(\by)\| &\geq \|\t{\by}^*-\proj_{\by_0+\spanop\mathcal{D}}(\by)\| - \|\by-\t{\by}^*\|, \\
        &\geq \|\t{\by}^*-\t{\by}_{\mathcal{D}}\| - \|\by-\t{\by}^*\|, \label{eq_poised_tmp3} \\
        &\geq \left(\alpha^*-\frac{1}{4}\right) \|\t{\by}-\t{\by}_{\mathcal{D}}\|, \\
        &> \left(\alpha^*-\frac{3}{4}\right) \|\t{\by}-\t{\by}_{\mathcal{D}}\|, \\
        &\geq \|\t{\by}^*-\t{\by}\| + \|\by-\t{\by}^*\|, \label{eq_poised_tmp5} \\
        &\geq \|\by-\t{\by}\|,
    \end{align}
    where \eqref{eq_poised_tmp3} follows as $\t{\by}^*$ is closer to $\t{\by}_{\mathcal{D}}$, its projection onto $\by_0+\spanop\mathcal{D}$, than $\proj_{\by_0+\spanop\mathcal{D}}(\by)$, and \eqref{eq_poised_tmp5} follows from \eqref{eq_poised_tmp4a}.
    Thus, just like $\t{\by}^*$, the point $\by$ also satisfies $\proj_C(\by) \notin \by_0+\spanop\mathcal{D}$.
    Since the directions $\bd_k$ in \algref{alg_dense_sampling} are dense in the unit sphere, there must be infinitely many $k$ for which $\by_0+r\bd_k \in B(\t{\by}^*,\|\t{\by}-\t{\by}_{\mathcal{D}}\|/4)$, and we are done.
\end{proof}


\subsubsection{Geometry Improving Algorithm} \label{subsec_geometry}
The last component we need is a geometry improving algorithm: that is, given an interpolation set $\{\by_0,\ldots,\by_n\}$ contained in $C\cap B(\by_0,r)$ for which $\{\by_t-\by_0 : t=1,\ldots,n\}$ is linearly independent, but where \eqref{eq_new_poised} does not hold, modify the set so that \eqref{eq_new_poised} also holds.
Together with \algref{alg_dense_sampling}, this gives us a method for constructing a $\Lambda$-poised interpolation set.

Fortunately, the poisedness-improving mechanism for the unconstrained case \cite[Algorithm 6.3]{Conn2009} remains suitable in the constrained case.
First, suppose our initial interpolation set is $\mathcal{Y}\defeq\{\by_0,\ldots,\by_n\}\subset C$.
Then
\begin{align}
    \bmat{1 & 0 \\ -\bee & I} \underbrace{\bmat{1 & \by_0^T \\ \vdots & \vdots \\ 1 & \by_n^T}}_{=:M} = \bmat{1 & \by_0^T \\ 0 & (\by_1-\by_0)^T \\ \vdots & \vdots \\ 0 & (\by_n-\by_0)^T}, \label{eq_M_to_L}
\end{align}
where $\bee\defeq(1,\ldots,1)^T\in\R^n$.
Hence $\det M=\det L \neq 0$ since $L$ has linearly independent rows by assumption.

Now, if $h:\R^n\to\R$ is any linear function, it must be equal to its linear interpolant (in Lagrange form) from \cite[Lemma 3.5]{Conn2009}; that is,
\begin{align}
    h(\by) = \sum_{t=0}^{n} h(\by_t) \ell_t(\by).
\end{align}
By taking the family of functions $h_0(\by)=1$ and $h_j(\by)=y_j$ for $j=1,\ldots,n$ we get
\begin{align}
    \underbrace{\bmat{1 & \cdots & 1 \\ \by_0 & \cdots & \by_n}}_{=M^T} \bmat{\ell_0(\by) \\ \vdots \\ \ell_n(\by)} = \bmat{1 \\ \by},
\end{align}
which is \cite[Eq.~(3.4)]{Conn2009} specialized to linear interpolation.
Hence by Cramer's rule, $|\ell_t(\by)|$ is the relative change in volume of the simplex given by $\mathcal{Y}$ when $\by_t$ is replaced by $\by$ (see \cite[p.~41]{Conn2009} for details).
The full poisedness improving algorithm is given in \algref{alg_poisedness_improving}, a simple extension of \cite[Algorithm 6.3]{Conn2009}.

\begin{algorithm}[tb]
\begin{algorithmic}[1]
\Require Constraint set $C\subset\R^n$ with projection operator $\proj_C:\R^n\to C$, starting point $\by_0\in C$, poisedness constant $\Lambda>1$, trust-region radius $\Delta$ and scaling factor $\beta\geq 1$.
\vspace{0.2em}
\State If a suitable initial interpolation set is not available, using \algref{alg_dense_sampling}, construct a set $\mathcal{Y}\subset C \cap B(\by_0,\beta\min(\Delta,1))$ with $|\det M| > 0$.
\For{$k=1,2,\ldots$}
    \State Construct the Lagrange polynomials $\ell_t$ for $t=0,\ldots,n$ for $\mathcal{Y}$ by solving the system \eqref{eq_lagrange_system} and calculate
    \begin{align}
        \Lambda_t := \max_{\by\in C\cap B(\by_0,\min(\Delta,1))} |\ell_t(\by)|. \label{eq_geometry_improving}
    \end{align}
    \If{$\Lambda_t > \Lambda$ for some $t$}
        \State For some $t_k$ and $\by_k$ with $|\ell_{t_k}(\by_k)| > \Lambda$, set $\mathcal{Y} = \mathcal{Y}\setminus\{\by_{t_k}\}\cup\{\by_k\}$.
    \Else
        \State \textbf{return} $\mathcal{Y}$.
    \EndIf
\EndFor
\end{algorithmic}
\caption{Method for constructing $\Lambda$-poised sets.}
\label{alg_poisedness_improving}
\end{algorithm}

\begin{theorem} \label{thm_poisedness_improving}
  If $C$ satisfies \assref{ass_feasible_set}, $\by_0\in C$, and we run \algref{alg_poisedness_improving} with $\Lambda>1$, $\Delta>0$ and $\beta\geq 1$, then it terminates in finite time, returning a set $\mathcal{Y} \subset C\cap B(\by_0,\beta\min(\Delta,1))$ which is $\Lambda$-poised for linear interpolation in $B(\by_0,\Delta)\cap C$.
\end{theorem}
\begin{proof}
  The initial $\mathcal{Y}$ is contained in $C\cap B(\by_0,\beta\min(\Delta,1))$ by assumption or by \thmref{thm_dense_sampling}.
  Any points added during the main loop of \algref{alg_poisedness_improving} must be in $C\cap B(\by_0,\min(\Delta,1))$, and so since $\beta\geq 1$ the output $\mathcal{Y}$ is also contained in $C\cap B(\by_0,\beta\min(\Delta,1))$.
  
  By construction, the output of \algref{alg_poisedness_improving} (if it terminates) must be $\Lambda$-poised for linear interpolation in $B(\by_0,\Delta)\cap C$.
  It remains to show that \algref{alg_poisedness_improving} terminates in finite time.
  
  Let $V_k$ be the volume of the simplex $\operatorname{vol}(\mathcal{Y})>0$ at the start of iteration $k$ of \algref{alg_poisedness_improving}.
  In particular, we have $V_1>0$ by assumption or by \thmref{thm_dense_sampling} (from $\det L \neq 0$).
  At each iteration of \algref{alg_poisedness_improving}, we modify $\mathcal{Y}$ in such a way that $V_{k+1}\geq \Lambda V_k$ by the reasoning above.
  So, if $V_{\max}$ is the maximum volume of all simplices contained in $C\cap B(\by_0,\min(\Delta,1))$, then \algref{alg_poisedness_improving} must terminate after at most $\lceil \log_{\Lambda}(V_{\max}/V_1)\rceil$ iterations.
\end{proof}

Of course, \thmref{thm_poisedness_improving} with \thmref{thm_fl_new} gives a way of constructing fully linear models in \algref{alg_cdfotr}.
We note that, compared to the unconstrained result \cite[Theorem 6.3]{Conn2009}, \thmref{thm_poisedness_improving} is a weaker result in that the number of iterations required for \algref{alg_poisedness_improving} to terminate is not uniformly bounded.

\subsection{Application to Composite Minimization} \label{subsec_comp}
In the previous sections, we outline how to construct linear interpolation models which satisfy the fully linear condition \defref{def_fl}.
Model-based DFO methods based on linear interpolation is particularly useful (e.g.~\cite{DFOGN}) for composite minimization, where the objective has the form
\begin{align}
    f(\bx) = F(\br(\bx)), \label{eq_f_composite}
\end{align}
where $\br:\R^n\to\R^N$ is a black-box function for which derivatives are unavailable, and $F:\R^N\to\R$ is some known function (e.g.~$F(\br)=\frac{1}{2}\|\br\|^2$ gives nonlinear least-squares minimization).

\begin{assumption} \label{ass_composite_smoothness}
  The composite objective \eqref{eq_f_composite} satisfies:
  \begin{enumerate}[label=(\alph*)]
      \item On the set $\cup_k B(\bx_k,\Delta_{\max})$, $\br$ is continuously differentiable and its Jacobian $J:\R^n\to\R^{N\times n}$ is Lipschitz continuous with constant $L_J$ and uniformly bounded by $J_{\max}$; and
      \item On the image under $\br$ of the set $\cup_k B(\bx_k,\Delta_{\max})$, $F$ is twice continuously differentiable with its gradient bounded by $G_{\max}$ and its Hessian bounded by $H_{\max}$, and is bounded below by $F_{\text{low}}$.
  \end{enumerate}
\end{assumption}

We note that under \assref{ass_composite_smoothness} that each component $r_i$ of $\br$ has $L_J$-Lipschitz continuous gradients, since
\begin{align}
    \|\grad r_i(\bx_1) - \grad r_i(\bx_2)\| = \|(J(\bx_1)-J(\bx_2))^T \bee_i\| \leq L_J\, \|\bx_1-\bx_2\|\, \|\bee_i\|.
\end{align}
Also, we have that $f$ is differentiable with $\grad f(\bx) = J(\bx)^T \grad F(\br(\bx))$, and $\grad f$ is Lipschitz continuous, from
\begin{align}
    \|\grad f(\bx_1) - \grad f(\bx_2)\| &\leq \|J(\bx_1)\| \, \|\grad F(\br(\bx_1)) - \grad F(\br(\bx_2))\| \nonumber \\
    &\qquad\qquad + \|J(\bx_1)-J(\bx_2)\|\,\|\grad F(\br(\bx_2))\|, \\
    &\leq J_{\max} H_{\max} \|\bx_1-\bx_2\| + L_J G_{\max}\|\bx_1-\bx_2\|.
\end{align}
All together, this means that our composite $f(\bx)$ \eqref{eq_f_composite} satisfies \assref{ass_smoothness} with $L_{\grad f} := J_{\max} H_{\max} + L_J G_{\max}$ and so the convergence theory from \secref{sec_wcc} applies.

\paragraph{Fully Linear Model Construction}
In this situation, we assume that we have the ability to evaluate $\br(\bx)$, and so can construct a linear interpolation model
\begin{align}
    \br(\by) \approx \bem(\by) := \b{c} + J (\by-\by_0),
\end{align}
by requiring that $\br(\by_t)=\bem(\by_t)$ for all $t=0,\ldots,n$ for some interpolation set $\{\by_0,\ldots,\by_n\}$.
This gives us the linear system (c.f.~\eqref{eq_interp_system})
\begin{align}
    M \bmat{\b{c}^T \\ J^T} = \bmat{\br(\by_0) \\ \vdots \\ \br(\by_n)}.
\end{align}
We then build a local quadratic model for $f(\bx)$ by taking a  Taylor series for $F$, evaluated at $\bem$:
\begin{align}
    f(\by) \approx m(\by) := c + \bg^T (\by-\by_0) + \frac{1}{2}(\by-\by_0)^T H (\by-\by_0), \label{eq_composite_model}
\end{align}
where the model gradient and Hessian are $\bg \defeq J^T \grad F(\b{c})$ and $H\defeq J^T \grad^2 F(\b{c}) J$.
The main result of this section is the following, which says that $m(\by)$ is a quadratic fully linear model for $f$, provided we can build linear models for $\br$ with interpolation over $\Lambda$-poised sets.

\begin{theorem} \label{thm_fl_new_composite}
  Suppose \assref{ass_composite_smoothness} holds and $\Delta\leq \Delta_{\max}$. 
  If the model $\bem$ is constructed from an interpolation set $\{\by_0,\ldots,\by_n\}$ which is $\Lambda$-poised for linear interpolation in $B(\by_0,\Delta)\cap C$ and $\|\by_t-\by_0\| \leq \beta \min(\Delta, 1)$, then $m(\by)$ \eqref{eq_composite_model} is a fully linear model for $f$ in $B(\by_0,\Delta)$, with
  \begin{align}
      \kappaef &= \frac{L_{\grad f}}{2} + G_{\max} \sqrt{N} \frac{n\Lambda \revision{L_J} \beta^2}{2} + H_{\max}\left(\sqrt{N} \Delta_{\max} \frac{n\Lambda \revision{L_J} \beta^2}{2} + J_{\max}\right)^2, \qquad \text{and} \\
      \kappaeg &=  G_{\max} \sqrt{N} \frac{n\Lambda \revision{L_J} \beta^2}{2}.
  \end{align}
\end{theorem}
\begin{proof}
  The interpolation condition $\br(\by_0)=\bem(\by_0)$ gives $\b{c}=\br(\by_0)$ and $c=f(\by_0)$.
  From \thmref{thm_fl_new} we have that each component model $m_t(\by)$ is fully linear for $r_t(\by)$ in $B(\by_0,\Delta)$.
  
  Now, fix $\bd\in C\cap B(\by_0,1)$.
  Then 
  \begin{align}
      |(\grad f(\by_0)-\bg)^T \bd| &= \left|(J(\by_0)^T \grad F(\br(\by_0)) - J^T \grad F(\b{c}))^T \bd\right|, \\
      &\leq \|\grad F(\br(\by_0))\| \, \|(J(\by_0)-J)\bd\|, \\
      &\leq G_{\max} \sqrt{N} \|(J(\by_0)-J)\bd\|_{\infty}, \\
      &\leq G_{\max} \sqrt{N} \frac{n\Lambda \revision{L_J} \beta^2}{2} \Delta, \label{eq_composite_tmp1}
  \end{align}
  where the last line follows from \eqref{eq_fl_g} applied to the component models $m_t$.
  Next, we instead fix $\bs\in C\cap B(\by_0,\Delta)$ and compute
  \begin{align}
      |f(\by_0+\bs) - m(\by_0+\bs)| &= \left|f(\by_0+\bs) - f(\by_0) - \bg^T \bs - \frac{1}{2}\bs^T H\bs\right|, \\
      &\leq \left|f(\by_0+\bs) - f(\by_0) - \grad f(\by_0)^T \bs\right| + \left|(\grad f(\by_0)-\bg)^T\bs\right| \nonumber \\
      &\qquad\qquad + \frac{1}{2}|\bs^T J^T \grad F^2(\br(\by_0)) J \bs|, \\
      &\leq \frac{L_{\grad f}}{2}\Delta^2 + \left|(\grad f(\by_0)-\bg)^T\bs\right| + H_{\max} \|J\bs\|^2. \label{eq_composite_tmp2}
  \end{align}
  Now, if $\Delta\geq 1$, we let $\bd=\Delta^{-1}\bs$ and so $\|\bd\|\leq 1$.
  Also, since $\|\bd\| \leq \|\bs\|$ and $\by_0,\by_0+\bs\in C$ we have $\by_0+\bd\in C$, which gives us
  \begin{align}
      \left|(\grad f(\by_0)-\bg)^T\bs\right| &= \Delta \left|(\grad f(\by_0)-\bg)^T\bd\right| \leq G_{\max} \sqrt{N} \frac{n\Lambda \revision{L_J} \beta^2}{2} \Delta^2, \label{eq_composite_tmp3}
  \end{align}
  where the last inequality follows from \eqref{eq_composite_tmp1}, and also
  \begin{align}
      \|J\bs\| &\leq \|(J-J(\by_0))\bs\| + \|J(\by_0)\bs\|, \\
      &\leq \sqrt{N} \Delta \|(J-J(\by_0))\bd\|_{\infty} + J_{\max} \Delta, \\
      &\leq \sqrt{N} \Delta^2 \frac{n\Lambda \revision{L_J} \beta^2}{2} + J_{\max}\Delta, \label{eq_composite_tmp4a} \\
      &\leq \left(\sqrt{N} \Delta_{\max} \frac{n\Lambda \revision{L_J} \beta^2}{2} + J_{\max}\right)\Delta, \label{eq_composite_tmp5}
  \end{align}
  where \eqref{eq_composite_tmp4a} follows from \eqref{eq_fl_g} applied to the component models $m_t$
  Instead, if $\Delta<1$, from \lemref{lem_fl_intermediate} we have
  \begin{align}
      \left|(\grad f(\by_0)-\bg)^T\bs\right| &\leq G_{\max} \sqrt{N} \|(J(\by_0)-J)\bs\|_{\infty} \leq  G_{\max} \sqrt{N} \frac{n \Lambda \revision{L_J} \beta^2}{2}\Delta^2, \label{eq_composite_tmp4}
  \end{align}
  and
  \begin{align}
      \|J\bs\| &\leq \|(J-J(\by_0))\bs\| + \|J(\by_0)\bs\|, \\
      &\leq \sqrt{N} \frac{n\Lambda \revision{L_J} \beta^2}{2} \Delta^2 + J_{\max} \Delta, \\
      &\leq \left(\sqrt{N} \Delta_{\max} \frac{n\Lambda \revision{L_J} \beta^2}{2} + J_{\max}\right)\Delta. \label{eq_composite_tmp6}
  \end{align}
  The result then follows from \eqref{eq_composite_tmp2} combined with either \eqref{eq_composite_tmp3} and \eqref{eq_composite_tmp5}, or \eqref{eq_composite_tmp4} and \eqref{eq_composite_tmp6}.
\end{proof}

We also note that the above proof also gives us an effective bound on the model Hessian, in the sense that
\begin{align}
    \|\bs^T H\bs\| \leq H_{\max}\left(\sqrt{N} \Delta_{\max} \frac{n\Lambda \revision{L_J} \beta^2}{2} + J_{\max}\right)^2 \Delta^2, \label{eq_composite_hess_bound}
\end{align}
for all $\bs\in C\cap B(\by_0,\Delta)$, whenever the interpolation set is $\Lambda$-poised.

\section{Numerical Results}
\subsection{Implementation}
To study the numerical performance of CDFO-TR with composite linear models (specifically for nonlinear least-squares problems), we extend the software package DFO-LS \cite{DFOLS}
to a new implementation we will call CDFO-LS\footnote{Version 1.3, available from \url{https://github.com/numericalalgorithmsgroup/dfols}}. 
CDFO-LS differs from DFO-LS primarily in its calculation of the trust-region step \eqref{eq_trs} in the presence of the constraint set $C$.
It also requires alternative methods for calculating  geometry-improving steps \eqref{eq_geometry_improving}, and the criticality measure \eqref{eq_pik}, however these may be viewed as special cases of \eqref{eq_trs} with a linear objective.
In all cases, this problem is convex (as the model Hessian in \eqref{eq_composite_model} has the form $J^T J$).\footnote{For geometry-improving steps \eqref{eq_geometry_improving}, we separately maximize $\pm \ell_t(\by)$.}
To efficiently solve these problems, we use an implementation of FISTA \cite[Chapter 10.7]{Beck2017} as given in \algref{alg_mfista}.
We terminate FISTA when either $\norm{\tilde{\bx}_{j+1} - \tilde{\bx}_j} \leq 10^{-12}$ or after $100n^2$ iterations, whichever is achieved first.
Projection onto the intersection of $C$ and a closed ball was performed using Dykstra's algorithm \cite{Dykstra1983,Dykstra1986} using the stopping criterion proposed in \cite{Birgin2005} with tolerance $10^{-10}$.

The other key difference between CDFO-LS and DFO-LS is the construction of an initial interpolation set. 
Our implementation uses \algref{alg_dense_sampling}, but where the sequence $\bd_k$ is constructed by first trying $\pm \bee_t$ for $t=1,\ldots,n$, and then generating the vectors $\bd_k$ i.i.d.~from a uniform distribution on the unit sphere.
DFO-LS only uses the choices $\pm \bee_t$, since it only allows the use of bound constraints.



\begin{algorithm}[tb]
\begin{algorithmic}[1]
\Require Constraint set $C\subset\R^n$, starting point $\bx_0\in C$, model gradient $\bg_0$ and Hessian $H$, and trust-region radius $\Delta$.
\State Set $\by_0 = \bx_0$, $L=\norm{H}$ if $H\neq 0$ (or $L=1$ otherwise), and  $t_0=1$.
\vspace{0.2em}
\For{$j=0,1,2,\ldots$}
    \State Project the gradient step onto $D := C \cap B(\bx_0,\Delta)$ and set $ \bx_{j+1} = \proj_D\left(\by_j - \frac{1}{L}\bg_j\right)$
    \If{the stopping conditions are satisfied}
        \vspace{0.5em}
        \State \textbf{return} the step $(\bx_{j+1} - \bx_j)$
        \vspace{0.5em}
    \Else
        \vspace{0.5em}
        \State Set $t_{j+1} = (1+\sqrt{1+4t_j^2})/2$ and $\by_{j+1} = \bx_{j+1} + \left(\frac{t_{j-1}}{t_{j+1}}\right)\left(\bx_{j+1} - \bx_j\right)$.
        \State Calculate the gradient at our new iterate $\bg_{j+1} = \bg_j + H(\by_{j+1} - \by_j)$. 
    \EndIf
\EndFor
\end{algorithmic}
\caption{Modified FISTA for solving \eqref{eq_trs}}
\label{alg_mfista}
\end{algorithm}


\paragraph{Implementation for box constraints}
The original implementation of DFO-LS can handle box constraints via a tailored trust-region subproblem solver \cite{Powell2009}.
For such problems we tested CDFO-LS (using projections) against the original DFO-LS, and found that the original implementation was slightly better.
Hence when only box constraints are present, CDFO-LS uses the original DFO-LS implementation.

\subsection{Solvers tested}
The numerical performance of CDFO-LS was measured against the following solvers:
\begin{itemize}
    \item COBYLA \cite{Powell1994}, a DFO solver for handling constrained general objectives based on a linear interpolation model, incorporated in the SciPy package; and
    \item pyNOMAD, the Python interface to NOMAD \cite{nomad2011}, an instance of Mesh Adaptive Direct Search (MADS) for black-box optimization of constrained general objectives.
\end{itemize}


Our tests on CDFO-LS used the default parameters. For all solvers, we used an initial trust-region radius of $\rho_0 = \Delta_0 = 0.1\max(\norm{\bm{x}_0}_{\infty},1)$. In performing tests with NOMAD, constraints were introduced as \textit{extreme barrier} constraints to prevent constraint violations. 
As NOMAD requires a feasible initial point, in this case we used $\proj_C(\bx_0)$. 


\subsection{Test problems and methodology}
The three solvers were tested on 5 problems\footnote{The names of the five problems chosen from \cite{Mor1981} are \textit{Biggs Exp6}, \textit{Dixon}, \textit{Gulf research and development}, \textit{Powell badly scaled}, and \textit{Wood}.} chosen out of the test suite from Moré, Garbow, and Hillstrom \cite{Mor1981}, and the entire test suite of 53 problems from Moré and Wild \cite{Mor2009}. All objectives were nonlinear least-squares functions with dimension $2 \leq n \leq 12$ and $2 \leq m \leq 65$.
For each problem, there were four types of constraints that were individually tested:
\begin{itemize}
    \item Unconstrained: $\bx\in\R^n$;
    \item Box: $\bx$ must satisfy (element-wise) $\frac{1}{10} \leq \bx \leq 20$;
    \item Ball: $\bx \in B(\bx_c, 6.9)$, $\bx_c = (5,5,\hdots,5)\in\R^n$; and
    \item Halfspace: we require $\bx \in \{\by : \bm{1}^Ty \leq 1\}$;
\end{itemize}
where $\bx \in \R^n$ is some point at which we evaluate the objective, and $\bm{1}$ is the vector of all ones. This means that in total there were $4 \times (53+5) = 232$ problems tested.


We additionally performed tests where stochastic noise was introduced when evaluating the residuals $r_i$. The two noise models implemented were as follows:
\begin{itemize}
    \item Multiplicative Gaussian noise: we evaluate residual $\tilde{r}_i(\bx) = r_i(\bx)(1+\epsilon)$; and
    \item Additive Gaussian noise: we evaluate residual $\tilde{r}_i(\bx) = r_i(\bx) + \epsilon$;
\end{itemize}
where $\epsilon \sim N(0,\sigma^2)$ i.i.d.~for each $i$ and $\bx$.

Prior to making comparisons, for each solver $\So$, every problem $p$, and an accuracy level $\tau \in (0,1)$, we determined the number of function evaluations $N_p(\So;\tau)$ required for a problem to be considered solved:
\begin{equation} \label{eq_evals_before_solved}
N_p(\So;\tau) := \text{\# objective evals before }\ f(\bx_k) \leq \E\left[f^* + \tau\left(f(\bx_0)-f^*\right)\right],
\end{equation}
where $f^*$ is the approximate true minimum $f(\bx^*)$ of the corresponding problem. If the problem is unconstrained, we source the true minimum from \cite{Mor1981} and \cite{Mor2009}. For all other cases, $f^*$ is taken to be the smallest objective value obtained by any of the three solvers. An important consideration in our comparisons is that the COBYLA algorithm allows constraints to be violated in search of a local minimizer. To keep our comparisons as fair as possible, any iteration for which COBYLA evaluated an infeasible point $\bx_k$ we recorded as $f(\bx_k) = \infty$. 
Moreover, if any solver is not able to obtain the true minimum for a corresponding $p$ and $\tau$ in the maximum allowed computational budget, then we set $N_p(\So,\tau)=\infty$.

To compare solvers, we use \textit{performance profiles} \cite{Dolan2002}, where we plot the number of objective evaluations $N_p(\So;\tau)$ normalized by the minimum objective evaluations needed by any solver $N^*_p(\tau)$:
\begin{equation} \label{eq_pi_performance}
\pi_{\So,\tau}(\alpha) := \frac{\abs{\{p \in \Ps : N_p(\So;\tau) \leq \alpha N_p^*(\tau)\}}}{\abs{\Ps}},\ \ \ \alpha \geq 1.
\end{equation}

\subsection{Test results}
In all tests, we generated performance profiles for accuracy levels $\tau=10^{-1}, \tau=10^{-3}$ and $\tau=10^{-5}$; in our noise models we used noise level $\sigma = 10^{-2}$. When testing CDFO-LS in the noisy case, we specified through a flag that the objective was noisy. With the flag set, CDFO-LS uses a different value of $\gammadec$ and changes the values of other parameters to better suit a noisy objective (see \cite{DFOLS} for details).

\begin{figure}[t]
  \centering
  \begin{subfigure}[b]{0.30\textwidth}
    \includegraphics[width=\textwidth]{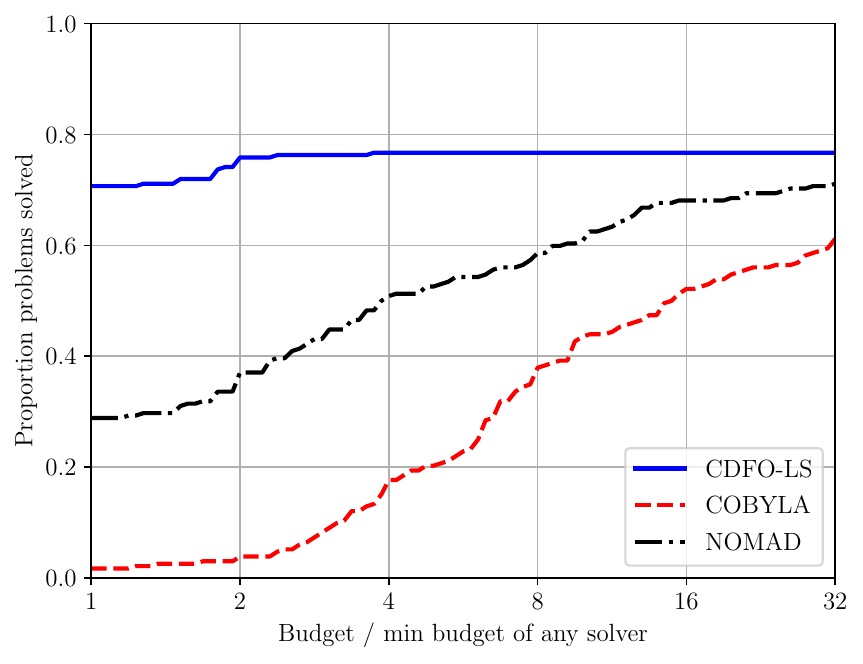}
    \caption{$\tau=10^{-1}$}
  \end{subfigure}
  ~
  \begin{subfigure}[b]{0.30\textwidth}
    \includegraphics[width=\textwidth]{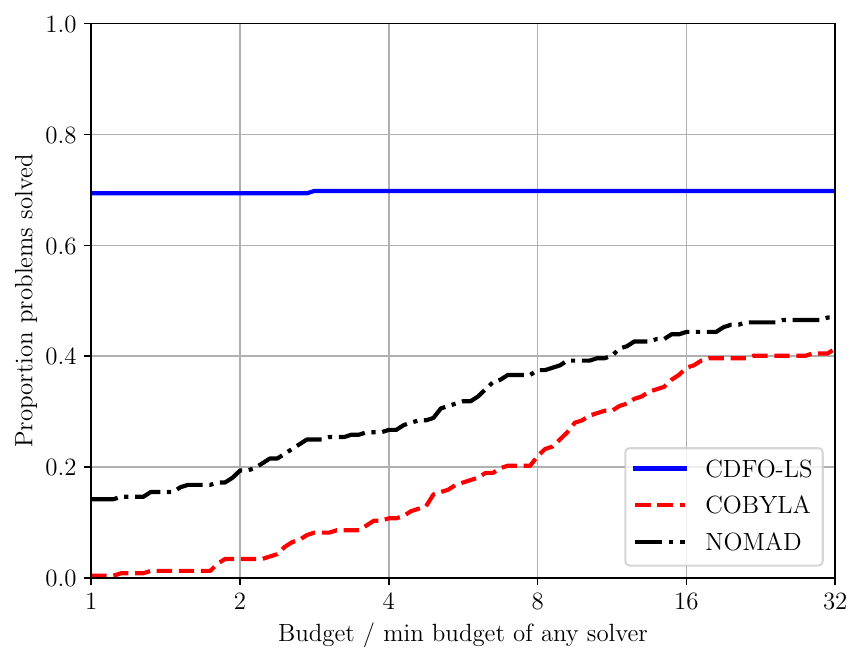}
    \caption{$\tau=10^{-3}$}
  \end{subfigure}
  ~
  \begin{subfigure}[b]{0.30\textwidth}
    \includegraphics[width=\textwidth]{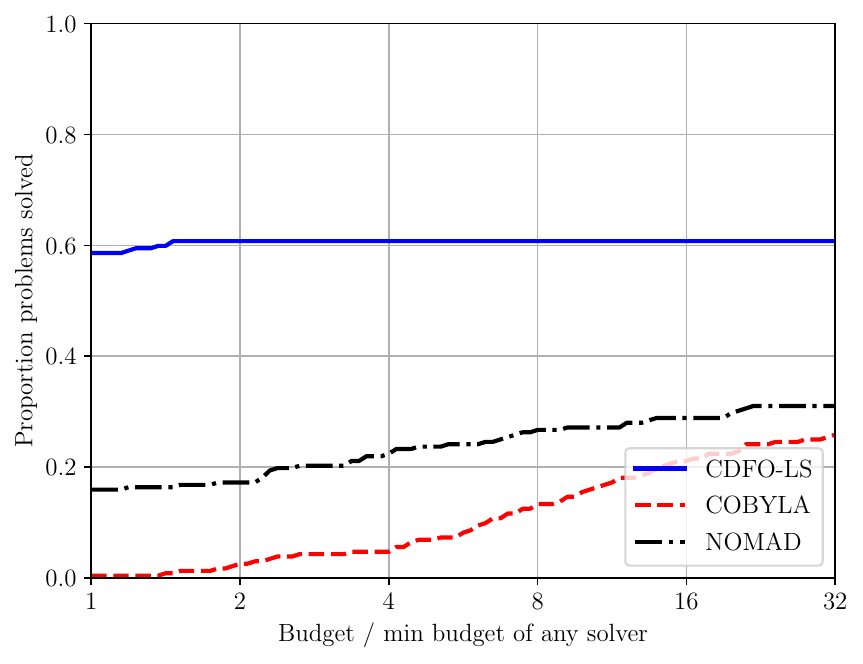}
    \caption{$\tau=10^{-5}$}
  \end{subfigure}
  \caption{Performance profiles comparing CDFO-LS against NOMAD and COBYLA on a collection of noiseless, constrained nonlinear least-squares problems.}
  \label{fig_noiseless_perf}
\end{figure}

\figref{fig_noiseless_perf} shows the performance profiles in the case that the objective is noiseless. We see in all accuracy settings, that CDFO-LS is able to solve a much greater proportion of problems in comparison to both COBYLA and NOMAD in a limited number of evaluations. 
This is to be expected, as it is able to explicitly exploit the least-squares problem structure, and has more information about the constraint set than NOMAD (a projection operator rather than a simple feasible/infeasible flag).
At the low accuracy level, as the number of evaluations increase, NOMAD manages to solve a similar proportion of problems to CDFO-LS; COBYLA performs worse, but is not far behind NOMAD.

As the accuracy level increases, the proportion of problems solved decreases across all solvers. However, the gap between the proportion of problems solved by CDFO-LS and NOMAD increases in favour of CDFO-LS. NOMAD performs more similarly to COBYLA at higher accuracy. Neither NOMAD nor COBYLA are able to solve many problems in comparison to CDFO-LS when $\tau = 10^{-5}$.

\begin{figure}[t]
  \centering
  \begin{subfigure}[b]{0.30\textwidth}
    \includegraphics[width=\textwidth]{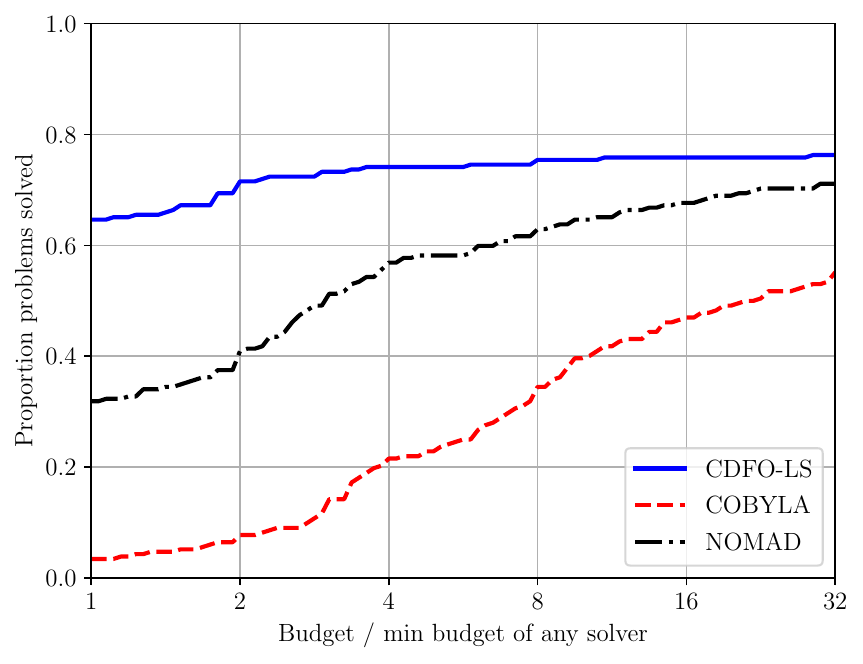}
    \caption{$\tau=10^{-1}$}
  \end{subfigure}
  ~
  \begin{subfigure}[b]{0.30\textwidth}
    \includegraphics[width=\textwidth]{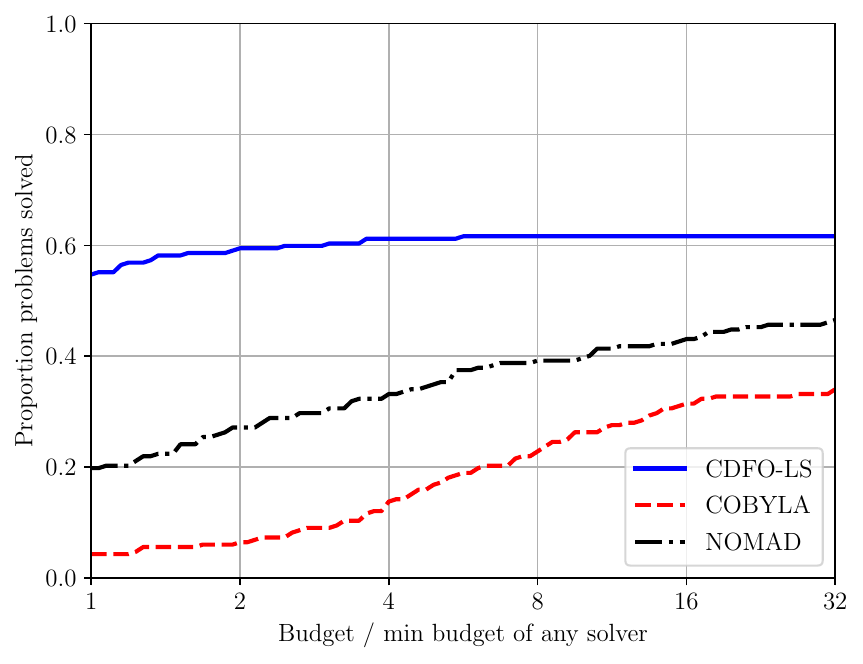}
    \caption{$\tau=10^{-3}$}
  \end{subfigure}
  ~
  \begin{subfigure}[b]{0.30\textwidth}
    \includegraphics[width=\textwidth]{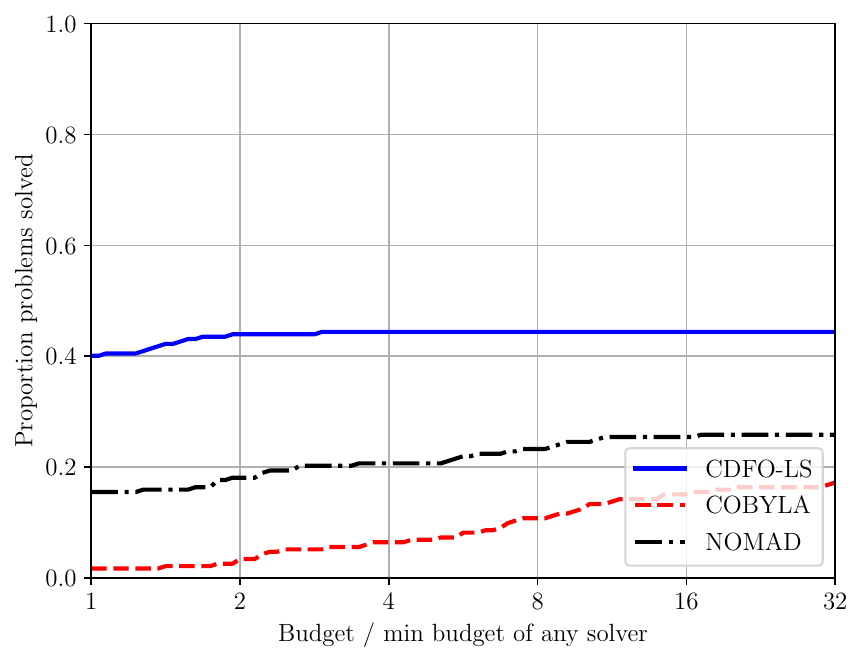}
    \caption{$\tau=10^{-5}$}
  \end{subfigure}
  \caption{Performance profiles comparing CDFO-LS against NOMAD and COBYLA on a collection of  constrained nonlinear least-squares problems in the presence of multiplicative noise.}
  \label{fig_mult_perf}
\end{figure}
\begin{figure}[t]
  \centering
  \begin{subfigure}[b]{0.30\textwidth}
    \includegraphics[width=\textwidth]{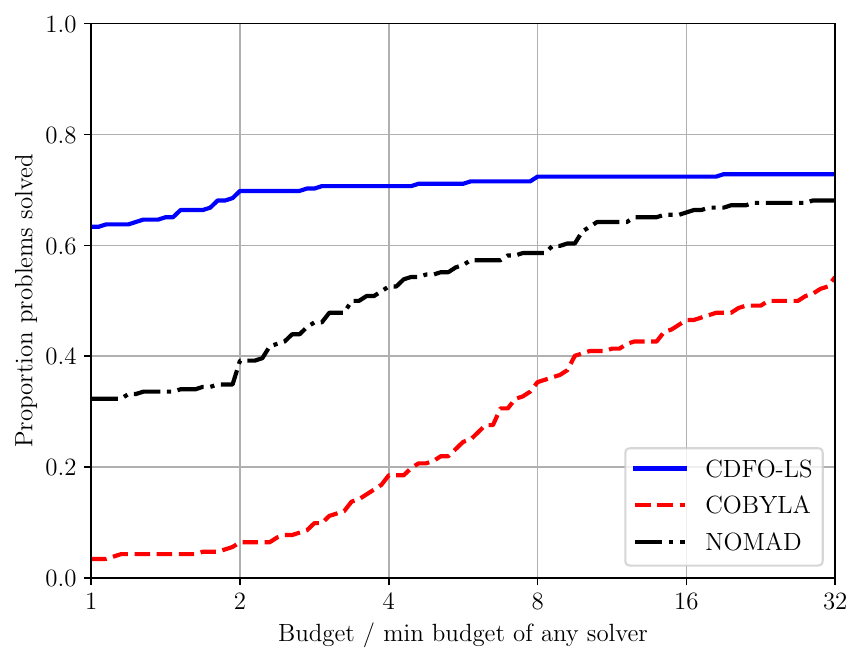}
    \caption{$\tau=10^{-1}$}
  \end{subfigure}
  ~
  \begin{subfigure}[b]{0.30\textwidth}
    \includegraphics[width=\textwidth]{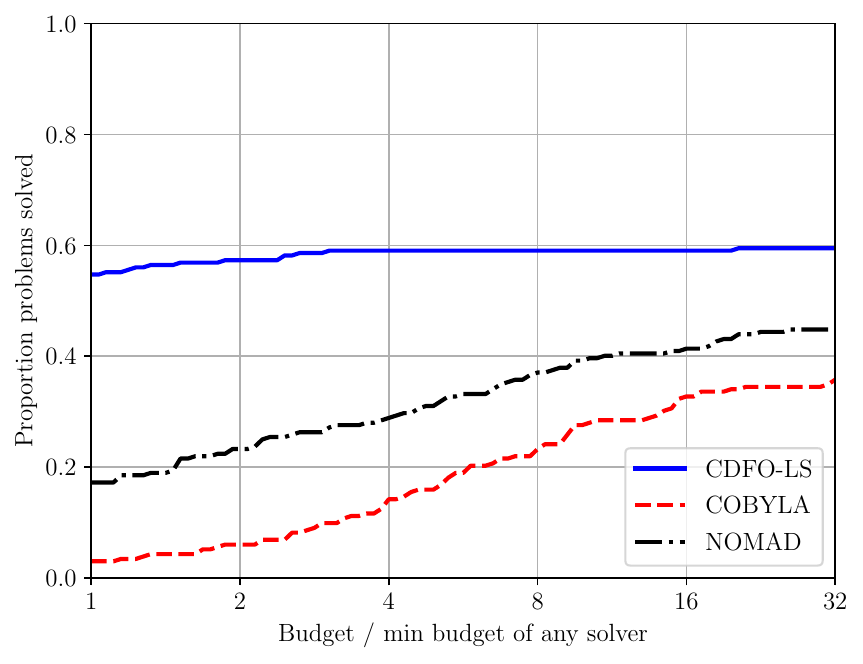}
    \caption{$\tau=10^{-3}$}
  \end{subfigure}
  ~
  \begin{subfigure}[b]{0.30\textwidth}
    \includegraphics[width=\textwidth]{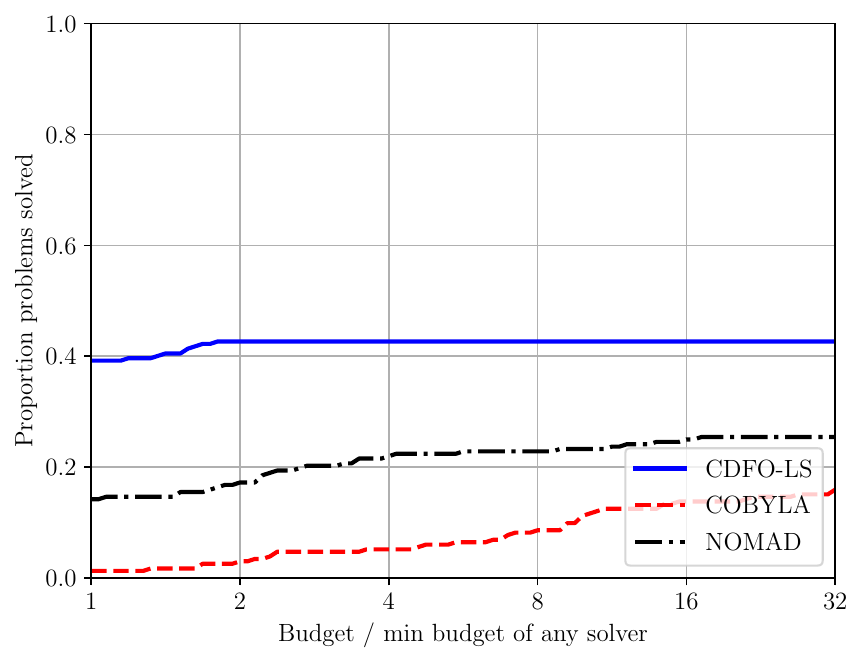}
    \caption{$\tau=10^{-5}$}
  \end{subfigure}
  \caption{Performance profiles comparing CDFO-LS against NOMAD and COBYLA on a collection of  constrained nonlinear least-squares problems in the presence of additive noise.}
  \label{fig_add_perf}
\end{figure}

In the case of noisy objectives, shown in Figures \ref{fig_mult_perf} and \ref{fig_add_perf}, the results are similar. However, in comparison to CDFO-LS, the performance of NOMAD and COBYLA is not as affected. Still, CDFO-LS significantly outperforms all other solvers, particularly at lower performance ratios $\alpha$.

\section{Conclusions and Future Work}
Current model-based DFO approaches for handling unrelaxable constraints consider only simple bound (\cite{Powell2009,Gratton2011} and \cite[Section 6.3]{Wild2009}) or linear inequality constraints \cite{Gumma2014}, or requires strong model accuracy assumptions that are difficult to satisfy using only feasible points \cite{Conejo2013}. 

We extend these approaches to handle a significantly broader class of constraints by assuming a convex constraint set accessed only via projections, while providing theoretical guarantees of convergence. 
A worst-case complexity analysis gives a bound of $\bigO(\epsilon^{-2})$ iterations to reach $\epsilon$ first-order optimality in the case of linear and composite linear interpolation models; the same result as in the unconstrained setting \cite{Garmanjani2016}. 
To achieve this, we introduce a weaker definition of a fully linear model (\defref{def_fl}) adapted from \cite{Conn2009a} for convex constraints. 
We also generalize the concept of $\Lambda$-poisedness  from \cite{Conn2007} to the constrained setting, which importantly allows $\Lambda$-poised interpolation sets to be constructed using only feasible points. 
We demonstrate how to verify and construct $\Lambda$-poised sets for linear and composite linear interpolation.

Numerically, we implement our method by extending  DFO-LS \cite{DFOLS} for nonlinear least-squares problems, a setting where linear interpolation models are sufficient to achieve have strong performance. Tests on low-dimensional test problems with a variety of convex constraints show that our implementation outperforms the general-purpose solvers COBYLA \cite{Powell1994} and NOMAD \cite{nomad2011}. 
We do note that neither solver is able to exploit the least-squares problem structure, and that COBYLA requires the functional form of the constraint sets, and NOMAD handles constraints via an extreme barrier approach. 

A natural extension to our work here is to develop convergence theory and complexity analysis for quadratic interpolation models, which would enable a practical implementation of our method for general nonconvex minimization. 

\bibliographystyle{siamplain}
\bibliography{references}
\end{document}